\RequirePackage{ifpdf}
\ifpdf % We are running pdfTeX in pdf mode
\documentclass[pdftex]{sigma}
\else
\documentclass{sigma}
\fi

\newcommand{\qbi}[2]{\left[
\begin{array}{c}
#1\\
#2
\end{array}
\right]}

\begin{document}

\allowdisplaybreaks

\renewcommand{\PaperNumber}{027}

\FirstPageHeading

\ShortArticleName{$q$-Wakimoto Modules and Integral Formulae}

\ArticleName{$\boldsymbol{q}$-Wakimoto Modules\\ and Integral Formulae of Solutions\\ of
the Quantum Knizhnik--Zamolodchikov Equations}

\Author{Kazunori KUROKI}

\AuthorNameForHeading{K. Kuroki}

\Address{Department of Mathematics, Kyushu University, Hakozaki 6-10-1, Fukuoka 812-8581, Japan}
\Email{\href{mailto:ma306012@math.kyushu-u.ac.jp}{ma306012@math.kyushu-u.ac.jp}}

\ArticleDates{Received October 31, 2008, in f\/inal form February 25,
2009; Published online March 07, 2009}

\Abstract{Matrix elements of intertwining operators between $q$-Wakimoto modules asso\-cia\-ted to the tensor product of representations of
 $U_q(\widehat{sl_2})$ with arbitrary spins are studied.
It is shown that they coincide with the Tarasov--Varchenko's formulae of the solutions of the qKZ equations.
The result generalizes that of the previous paper [Kuroki~K., Nakayashiki~A., \href{http://dx.doi.org/10.3842/SIGMA.2008.049}{{\it SIGMA} {\bf 4} (2008), 049, 13~pages}].}

\Keywords{free f\/ield; vertex operator; qKZ equation; $q$-Wakimoto module}

\Classification{81R50; 20G42; 17B69}

\section{Introduction}

 In \cite{nakayashiki} the integral formulae of the quantum Knizhnik--Zamolodchikov (qKZ) equations \cite{frenkel2}
for the tensor product of spin $1/2$ representation of $U_q({sl_2})$ arising from
$q$-Wakimoto modules have been studied.
The formulae are identif\/ied with those of Tarasov--Varchenko's formulae.
The aim of this paper is to generalize the results to the case of tensor product of representations with arbitrary spins.

It is known that certain matrix elements of intertwining operators between $q$-Wakimoto modu\-les satisfy the qKZ equation~\cite{frenkel2,matsuo}.
 Thus it is interesting to compute those matrix elements explicitly.
In \cite{jimbo} two kinds of intertwining operators were introduced, type I and type II.
They were def\/ined according as the position of evaluation representations.
In the application to the study of solvable lattice models two types of operators have their own roles.
Type I and type II operators correspond to states and particles respectively.
The properties of traces  exhibit very dif\/ferent structure.
However as far as the matrix elements are concerned they are not expected to be very dif\/ferent \cite{jimbo}.

In \cite{nakayashiki} a computation of matrix elements has been carried out in the case of type I ope\-ra\-tor and the tensor product of 2-dimensional vector representation of $U_q(sl_2)$
generalizing the result of~\cite{matsuo} (see the previous paper \cite{nakayashiki}).
In this paper we compute matrix elements for the composition of the type I
intertwining operators \cite{jimbo} associated to f\/inite dimensional irreducible representations
of $U_q(sl_2)$. We perform certain multidimensional integrals and sums explicitly. It is shown that the formulae thus obtained coincide with those of
Matsuo \cite{M0}, Tarasov and Varchenko \cite{tarasov2} without the term corresponding to the deformed cycles.

To obtain actual matrix elements of intertwining operators it is necessary to specify certain contours of
integration associated to screening operators.
We do not consider this problem in this paper.
To f\/ind integration contours describing each composition of intertwining operators
is an important open problem.
We also remark that the formulae for type II intertwining operators are not obtained
in this paper. The computation of them looks quite dif\/ferent from that for type I case as opposed to the expectation.
It is interesting to f\/ind the way to get a similar result for matrix elements in the case of type II
operators.

 The paper is organized in the following manner.
The construction of the solutions of the qKZ equations due to Tarasov and Varchenko is reviewed in Section~\ref{section2}.
In Section~\ref{section3} a free f\/ield construction of intertwining operator is reviewed.
The formulae for the matrix elements of some operators are calculated in Section~\ref{section4}.
The main theorem of this paper is stated in this section.
In Section~\ref{section5} the proof of the main theorem is given.
The evaluation representation of $U_q(\widehat{sl_2})$ is explicitly described in  Appendix~\ref{appendixA}.
Appendix~\ref{appendixB} gives the explicit form of the $R$-matrix in special cases.
The explicit forms of the operators which appear in Section~\ref{section3} are given in Appendix~\ref{appendixC}.
Appendix~\ref{appendixD} contains the list of OPE's which is necessary to derive the integral formulae.

\section[Tarasov-Varchenko's formulae]{Tarasov--Varchenko's formulae}\label{section2}

We review Tarasov--Varchenko's formula for solutions of  the qKZ equations. In this paper we assume that $q$ is a complex number such that $|q|<1$. We mainly follow the notation of~\cite{tarasov2}.
For a nonnegative integer $l$ let $V^{(l)}=\bigoplus_{i=0}^l {\mathbb{C}}v_i^{(l)}$ be the  $l+1$ dimensional irreducible $U_q(sl_2)$-module and $V^{(l)}_z=V^{(l)}\otimes {\mathbb{C}}[z, z^{-1}]$ the evaluation representation of $U_q(\widehat{sl_2})$ on $V^{(l)}$.
The action of $U_q(\widehat{sl_2})$ on $V^{(l)}_z$ is given in Appendix~\ref{appendixA}.
Let $l_1$ and $l_2$ be nonnegative integers and $R_{l_1, l_2}(z)\in {\rm End}(V^{(l_1)}\otimes V^{(l_2)})$ the trigonometric quantum $R$-matrix uniquely determined  by
the following conditions:
\begin{gather*}
(i)~ \ \ P R_{l_1, l_2}(z)  \ \text{commutes with} \ U_q(\widehat{sl_2}),
\\
(ii) \ \ P R_{l_1, l_2}(z) \big(v_0^{(l_1)}\otimes v_0^{(l_2)}\big)=v_0^{(l_2)}\otimes v_0^{(l_1)},
\end{gather*}
where $P : V^{(l_1)}\otimes V^{(l_2)}\to V^{(l_2)}\otimes V^{(l_1)}$ is a linear map given by
\begin{gather*}
P(v \otimes w)=w \otimes v.
\end{gather*}
The explicit form of the $R$-matrix is given in Appendix~\ref{appendixB} in case $l_1=1$ or $l_2=1$.
We set
\begin{gather*}
 \widehat{R}_{l_i, l_j}(z)=\rho_{l_i, l_j}(z)\widetilde{R}_{l_i, l_j}(z),
\qquad \widetilde{R}_{l_i, l_j}(z)=(C_{l_i}\otimes C_{l_j}){R}_{l_i, l_j}(z)(C_{l_i}\otimes C_{l_j}),
\\
\rho_{l_i, l_j}(z)=q^{\frac{l_i l_j}{2}}\frac{(q^{l_i+l_j+2}z^{-1}; q^4)_{\infty}(q^{-l_i-l_j+2}z^{-1}; q^4)_{\infty}}
{(q^{-l_i+l_j+2}z^{-1}; q^4)_{\infty}(q^{l_i-l_j+2}z^{-1}; q^4)_{\infty}},
\\
C_{l} v_{\epsilon}^{(l)}=v_{l-\epsilon}^{(l)}\qquad \big(v_\epsilon^{(l)}\in V^{(l)}\big),
\end{gather*}
where for a complex number $a$ with $|a|<1$
\begin{gather*}
(z; a)_{\infty}=\prod_{i=0}^{\infty}\big(1-a^i z\big).
\end{gather*}
Let $k$ be a complex number. We set
\begin{gather*}
p=q^{2(k+2)}.
\end{gather*}
We assume that $p$ satisf\/ies $|p|<1$.
Let $T_j$ denote the $p$-shift operator of $z_j$,
\begin{gather*}
T_j f(z_1, \dots, z_n)=f(z_1, \dots, p z_j, \dots z_n).
\end{gather*}
Let $l_1, \dots, l_n$ and $N$ be nonnegative integers.
The qKZ equation for a $V_{l_1}\otimes \cdots \otimes V_{l_n}$-valued function $\Psi(z_1, \dots, z_n)$
is
\begin{gather}
T_j\Psi=\widehat{R}_{j, j-1}(p z_j/z_{j-1})\cdots \widehat{R}_{j, 1}(p z_j/z_{1})\kappa^{\frac{h_j}{2}}\widehat{R}_{j, n}(z_j/z_{n})\cdots
\widehat{R}_{j, j+1}(z_j/z_{j+1})\Psi,
\label{qKZ}
\end{gather}
where $\kappa$ is a complex parameter,
$\widehat{R}_{i, j}(z)$ signif\/ies that $\widehat{R}_{l_i, l_j}(z)$ acts on the $i$-th and $j$-th components of the tensor product and
$\kappa^{h_j}$ acts on $j$-th component as
\begin{gather*}
\kappa^{\frac{h_j}{2}}v_{m}^{(l_j)}=\kappa^{\frac{l_j-2m}{2}} v_{m}^{(l_j)}.
\end{gather*}
We set
\begin{gather*}
(z)_{\infty}=(z; p), \qquad
\theta(z)=(z)_{\infty}\big(pz^{-1}\big)_{\infty}(p)_{\infty}.
\end{gather*}

Consider a sequence $(\nu)=(\nu_1, \dots, \nu_n)$ satisfying $0 \le\nu_i \le l_i$ for all $i$
and $N=\sum\limits_{i=1}^n \nu_i$.
Let $r=\sharp\{i\,|\,\nu_i\ne 0\}$, $\{i\,|\,\nu_i\ne 0\}=\{k(1)<\dots<k(r)\}$ and
$n_i=\nu_{k(i)}$. We set
\begin{gather*}
w_{(\nu)}(t, z)=\prod_{a<b}\frac{t_a-t_b}{q^{-2}t_a-t_b}
\sum_{\substack{\Gamma_1\sqcup\dots \sqcup \Gamma_r=\{1, \dots, N\} \\ |\Gamma_s|=n_s (s=1, \dots, r)}}
\left(\prod_{\substack{1\le i<j\le r  \\  a\in \Gamma_i, b\in \Gamma_j}}
\frac{q^{-2}t_{a}-t_{b}}{t_{a}-t_{b}}\right)
\nonumber
\\
\phantom{w_{(\nu)}(t, z)=}{} \times
\prod_{b\in \Gamma_s}\Bigg(
\frac{t_{b}}{t_{b}-q^{-l_{k(i)}}z_{k(i)}}\prod_{j<k(i)}\frac{q^{-l_{j}}t_{b}-z_j}{t_{b}-q^{-l_{j}}z_j}
\Bigg).
%\label{w}
\end{gather*}
The elliptic hypergeometric space ${\cal{F}}_{\rm ell}$ is the space of functions $W(t, z){=}W(t_1, \dots, t_N, z_1, \dots, z_n)$
of the form
\begin{gather*}
W=Y(z)\Theta(t, z)\frac{1}{\prod\limits_{j=1}^n\prod\limits_{a=1}^N \theta(q^{l_j} t_a/z_j)}
\prod_{1\le a<b \le N}\frac{\theta(t_a/t_b)}{\theta(q^{-2}t_a/t_b)}
\end{gather*}
satisfying the following conditions:
\begin{enumerate}\itemsep=0pt
\item[$(i)$] $Y(z)$ is meromorphic on $({\mathbb{C}}^{\ast})^n$ in $z_1, \dots, z_n$, where
${\mathbb{C}}^{\ast}={\mathbb{C}}\setminus \{0\}$;

\item[$(ii)$] $\Theta(t, z)$ is holomorphic on $({\mathbb{C}}^{\ast})^{n+N}$ in $t_1, \dots, z_n$ and symmetric in $t_1, \dots, t_N$;

\item[$(iii)$]
$T_a^tW/W=\kappa q^{-2N+4a-2} \prod\limits_{i=1}^n q^{l_i}$, $T_j^zW/W=q^{-l_j N}$,
where $T_a^tW=W(t_1, \dots, p t_a, \dots, t_N, z)$ and $T_j^zW=W(t, z_1, \dots, p z_j, \dots, z_n)$.
\end{enumerate}

Def\/ine the phase function $\Phi(t, z)$ by
\begin{gather*}
\Phi(t, z)=\left(\prod_{a=1}^{N}\prod_{i=1}^{n}\frac{(q^{l_i} t_a /z_i)_{\infty}}{(q^{-l_i} t_a /z_i)_{\infty}}\right)
\left(\prod_{a<b}\frac{(q^{-2}t_a/t_b)_{\infty}}{(q^{2}t_a/t_b)_{\infty}}\right).
\end{gather*}

For $W\in {\cal{F}}_{\rm ell}$ let
\begin{gather}
I(w_{(\epsilon)}, W)=\int_{\widetilde{\mathbb{T}}^N}\prod_{a=1}^N\frac{d t_a}{t_a}\Phi(t, z) w_{(\epsilon)}(t, z)
W(t, z),
\label{TV-1}
\end{gather}
where $\widetilde{\mathbb{T}}^N$ is a suitable deformation of the torus
\begin{gather*}
\mathbb{T}^N=\{(t_1, \dots, t_N)\,|\, |t_i|=1, 1\le i \le N\},
\end{gather*}
specif\/ied as follows.
The integrand has simple poles at
\begin{gather*}
 t_a/z_j=\big(p^s q^{-l_j}\big)^{\pm 1},\qquad  s\ge 0, \quad 1\le a \le N, \quad 1\le j \le n,
\\
 t_a/t_b=\big(p^s q^2\big)^{\pm 1},\qquad  s\ge 0, \quad 1\le a< b\le N.
\end{gather*}
The contour of integration in $t_a$ is a simple closed curve which rounds the origin in the counterclockwise direction and separates the following two sets
\begin{gather*}
 \big\{p^s q^{-l_j}z_j, p^s q^2 t_b|s\ge 0, 1\le j \le N, a<b\big\},
\\
\big\{p^{-s} q^{l_j}z_j, p^{-s} q^{-2} t_b|s\ge 0, 1\le j\le N, a<b\big\}.
\end{gather*}
Let $L$ be a complex number and
\begin{gather*}
\kappa=q^{-2\big(L+\sum\limits_{i=1}^n\frac{l_i}{2}-N+1\big)}.
\end{gather*}
Then
\begin{gather}
\Psi_W=\Bigg(\prod_{i=1}^n z_i^{a_i}\Bigg)
\Bigg(\prod_{i<j}\xi_{l_i, l_j}(z_i/z_j)\Bigg)
\sum_{(\epsilon)}I(w_{(-\epsilon)}, W) v^{(l_1)}_{\epsilon_1}\otimes \dots \otimes v^{(l_n)}_{\epsilon_n}
\label{TV-2}
\end{gather}
is a solution of the qKZ equation (\ref{qKZ}) for any $W \in {\cal{F}}_{\rm ell}$
where $(-\epsilon)=(l_1-\epsilon_1, \dots, l_n-\epsilon_n)$ and
\begin{gather*}
a_i={\frac{l_i}{2(k+2)} }\Bigg(L+\sum_{j=1}^n{l_j}-\frac{l_i}{2}-N+1\Bigg),
\\
\xi_{l_i, l_j}(z)
=\frac{\left(p q^{l_i+l_j+2}z^{-1} ; q^4 , p \right)_{\infty}
\left(p q^{-l_i-l_j+2}z^{-1} ; q^4 , p \right)_{\infty}}
{\left(p q^{l_i-l_j+2}z^{-1} ; q^4 , p \right)_{\infty}
\left(p q^{-l_i+l_j+2}z^{-1}  ; q^4 , p \right)_{\infty}},
\\
(z; p, q)=\prod_{i=0}^{\infty}\prod_{j=0}^{\infty}(1-p^i q^j z).
\end{gather*}

\section[Free field realizations]{Free f\/ield realizations}\label{section3}

We brief\/ly review the free f\/ield construction of the representation of the $U_q(\widehat{sl_2})$ of level $k$ \cite{abada,matsuo,shiraishi} and
intertwining operators \cite{BW,kato,konno}.
We mainly follow the notation of \cite{kato}.
We set
\begin{gather*}
[n]=\frac{q^n-q^{-n}}{q-q^{-1}}.
\end{gather*}
Let $k$ be a complex number and $\{a_n, b_n ,c_n, \tilde{a}_0, \tilde{b}_0, \tilde{c}_0, Q_a, Q_b, Q_c\,|\,n\in {\mathbb{Z}}_{\ge 0}\}$ satisfy
\begin{gather*}
 [a_n, a_m]=\delta_{m+n, 0}\frac{[(k+2)n][2n]}{n},\qquad [\tilde{a}_0, Q_a]=2(k+2),
\\
 [b_n, b_m]=\delta_{m+n, 0}\frac{-[2n]^2}{n},\qquad [\tilde{b}_0, Q_b]=-4,
\\
[c_n, c_m]=\delta_{m+n, 0}\frac{[2n]^2}{n},\qquad [\tilde{c}_0, Q_c]=4.
\end{gather*}
Other combinations of elements are supposed to commute. Set
\begin{gather*}
N_{\pm}={\mathbb{C}}[a_n, b_n, c_n \,|\, \pm n > 0].
\end{gather*}
Let $r$ be a complex number and $s$ an integer.
The Fock module $F_{r, s}$ is def\/ined to be the free~$N_-$ module of rank one generated by the vector $|r, s\rangle$ satisfying
\begin{gather*}
N_{+}|r, s\rangle=0,
\qquad
\tilde{a}_0|r,s\rangle=r|r, s\rangle,
\qquad
\tilde{b}_0|r,s\rangle=-2s|r, s\rangle,
\qquad
\tilde{c}_0|r,s\rangle=-2s|r, s\rangle.
\end{gather*}
We set
\begin{gather*}
F_r=\oplus_{s\in {\mathbb{Z}}}F_{r, s}.
\end{gather*}
The right Fock module $F_{r, s}^{\dagger}$ and $F_{r}^{\dagger}$ are similarly def\/ined using the vector
$\langle r, s|$ satisfying the conditions
\begin{gather*}
\langle  r,s|N_{-}=0,
\qquad
 \langle r,s|\tilde{a}_0=r\langle r,s|,
\qquad
 \langle r,s|\tilde{b}_0=-2s\langle r,s|,
\qquad
\langle r,s|\tilde{c}_0=-2s\langle r,s|.
\end{gather*}
Notice that $F_r$ and $F_{r}^{\dagger}$ have left and right $U_q(\widehat{sl_2})$-module structure respectively~\cite{matsuo, shiraishi}.

Let
\begin{gather*}
|L\rangle=|L, 0\rangle\in F_{L, 0},\qquad \langle L|=\langle L, 0|\in F_{L, 0}^{\dagger}.
\nonumber
\end{gather*}
They become left and right highest weight vectors of $U_q(\widehat{sl_2})$ with the weight $L\Lambda_1+(k-L)\Lambda_0$ respectively, where $\Lambda_0$ and $\Lambda_1$ are fundamental weights of $\widehat{sl_2}$.

We consider operators
\begin{gather*}
 \phi_m^{(l)}(z): \ F_{r, s}\to F_{r+l, s+l-m},
\qquad
 J^-(u): \ F_{r, s}\to F_{r, s+1},
\qquad
 S(t): \ F_{r, s}\to F_{r-2, s-1},
\end{gather*}
the explicit forms of which are given in Appendix~\ref{appendixC}. We set
\begin{gather*}
\phi_l^{(l)}(z)=\phi_l(z)
\end{gather*}
for simplicity.
The operator $\phi_m^{(l)}(z)$ is used to construct the vertex operator for $U_q(\widehat{sl_2})$:
\begin{gather*}
 \phi^{(l)}(z): \ W_r\to W_{r+l}{\otimes} V^{(l)}_z,
\qquad
 \phi^{(l)}(z)=\sum_{m=0}^l \phi^{(l)}_m(z)\otimes v_m^{(l)},
\end{gather*}
where $W_r$ is a certain submodule of $F_r$ called $q$-Wakimoto module \cite{matsuo}.

The operator $J^-(u)$ is a generating function of a part of generators of the Drinfeld realization for $U_q(\widehat{sl_2})$
at level $k$.

The operator $S(t)$ commutes with $U_q(\widehat{sl_2})$ modulo total dif\/ferences. Here modulo total dif\/ferences means
modulo functions of the form
\begin{gather*}
_{k+2}\partial_z f(z)=\frac{f(q^{k+2}z)-f(q^{-(k+2)}z)}{(q-q^{-1})z}.
\end{gather*}

Consider
\begin{gather*}
F(t,z)=\langle L+\sum_{i=1}^n l_i -2N|\phi^{(l_1)}(z_1)\cdots\phi^{(l_n)}(z_n) S(t_N)\cdots S(t_1)|L\rangle
%\label{F}
\end{gather*}
which is a function taking the value in $V^{(l_1)}\otimes\dots\otimes V^{(l_n)}$.
Let
\begin{gather*}
\triangle_j=\frac{j(j+2)}{4(k+2)}.
\end{gather*}
Set
\begin{gather*}
\widehat{F}=\left(\prod_{i=1}^n z_i^{\triangle_{L+\sum\limits_{j=i}^{n}l_j-2N}-\triangle_{L+\sum\limits_{j=i+1}^{n}l_j-2N}}\right)F
=\left(\prod_{i=1}^n z_i^{\frac{l_i}{2(k+2)} \big(L+\sum\limits_{i<j}l_j-2N+\frac{l_i+2}{2}\big)}\right)F.
%\label{Fhat}
\end{gather*}

 Then the function $\widehat{F}(t,z)$ satisf\/ies qKZ equation (\ref{qKZ}) with $\kappa=q^{-2\big(L+\frac{\sum\limits_{i=1}^n l_i}{2}-N+1\big)}$ modulo total dif\/ferences \cite{matsuo}.

\section{Integral formulae}\label{section4}

Def\/ine the components of $F(t,z)$ by
\begin{gather*}
F(t,z)=\sum_{\substack{\nu_i\in\{0, \dots, l_i\} \\ 1\le i \le n}}
F^{(\nu)}(t,z) v_{\nu_1}^{(l_1)}\otimes \dots \otimes v_{\nu_n}^{(l_n)},
\end{gather*}
where  $(\nu)=(\nu_1, \dots, \nu_n)$. By the conditions on weights  $F^{(\nu)}(t,z)=0$ unless
\begin{gather*}
\sum_{i=1}^n(l_i-\nu_i)=N
\end{gather*}
is satisf\/ied. We assume this condition once for all.
Let
\begin{gather*}
\sharp\{i\, |\, \nu_i\ne l_i\}=r,\qquad \{i\, |\, \nu_i\ne l_i\}=\{k(1)<\cdots<k(r)\},\\
n_i=l_{k(i)}-\nu_{k(i)} \qquad (1\le i \le r).
\end{gather*}

The main result of this paper is
\begin{theorem}\label{main}
We have
\begin{gather*}
F^{(\nu)}(t, z) = A^{(\nu)}(t, z)\left(\prod_{i=1}^n z_i^{\frac{l_i}{2(k+2)} \big(L-2N+\sum\limits_{i<j}l_j\big)}\right)
\Bigg(\prod_{i<j}\xi_{l_i, l_j}(z_i/z_j)\Bigg)
\Phi(t, z)w_{(-\nu)}(t, z)
,
\end{gather*}
where $(-\nu)=(l_{1}-\nu_1, \dots, l_{n}-\nu_n)$,
$n_i=l_{k(i)}-\nu_{k(i)}$ and
\begin{gather*}
A^{(\nu)}(t, z) = q^{-NL}q^{\frac{3N(N-1)}{2}-\big(\sum\limits_{i=1}^n l_i\big)N}
q^{\frac{1}{2(k+2)}\big(k{\sum\limits_{i<j}l_i l_j}+k(L-2N)\sum\limits_{i=1}^nl_i+4LN-4N(N-1)\big)}
\\
\phantom{A^{(\nu)}(t, z) =}{} \times
\left(\frac{1}{q-q^{-1}}\right)^{N}
\sum_{(\nu)}
\left\{\prod_{s=1}^{r}q^{\big(\sum\limits_{t=s+1}^r n_t\big) n_s-l_{k(s)}n_s} \right\}
\left\{\prod_{s=1}^r \prod_{i=0}^{n_s-1}\big(1-q^{2(l_{k(s)}-i)}\big)\right\}
\\
\phantom{A^{(\nu)}(t, z) =}{}\times
\left(\prod_{a=1}^N t_a^{\frac{2}{k+2}(a-1)-\frac{1}{k+2}L-1}\right).
\end{gather*}
\end{theorem}

The formula for $F^{(\nu)}(t, z)$ is of the form of (\ref{TV-1}), (\ref{TV-2}).  More precisely in Tarasov--Varchenko's formula (\ref{TV-1}), (\ref{TV-2}), $W$ can be written as
\begin{gather*}
W=
\left(\prod_{i=1}^{n}z_i^{\frac{l_i}{2(k+2)}\big(L-3N-\sum\limits_{j< i}l_j+\sum_{i<j}l_j\big)}\right)
\left(\prod_{a=1}^{N}t_a \right)
A^{(\nu)}(t, z)  W'
\end{gather*}
for suitable $W'$. This $W'$ specif\/ies an intertwiner. In this paper we don't consider the problem on specifying $W'$.

To prove Theorem~\ref{main} let us begin by writing down the formula obtained by the free f\/ield description of operators
$\phi_l(z)$, $J^{-}(u)$, $S(t)$ given in Appendix~\ref{appendixC}.
Let $(\epsilon)=(\epsilon_1, \dots, \epsilon_N), (\mu)=(\mu_{1, 1}, \dots, \mu_{1, n_1}, \dots, \mu_{r, n_r})\in \{0, 1\}^N$. Then $F^{(\nu)}(t, z)$ can be written as
\begin{gather*}
F^{(\nu)}(t, z) = (-1)^N \big(q-q^{-1}\big)^{-2N} \prod_{i=1}^{r}\frac{1}{[n_i]!} \prod_{a=1}^N t_a^{-1}
\\
\phantom{F^{(\nu)}(t, z) =}{}
\times
\sum_{\epsilon_i, \mu_{i_1, i_2}=\pm}\prod_{i=1}^N  \epsilon_i
\oint\left(\prod_{\genfrac{}{}{0pt}{}{1\le i_1 \le r}{1\le i_2\le n_{i_1}}}\mu_{i_1, i_2}
\frac{d u_{i_1, i_2}}{2 \pi i u_{i_1, i_2}}\right)F_{(\epsilon)(\mu)}^{(\nu)}(t, z| u),
\end{gather*}
where
\begin{gather*}
F_{(\epsilon)(\mu)}^{(\nu)}(t, z| u)
=
\Big\langle L+\sum_{i=1}^{n}l_i-2N|\phi_{l_1}(z_1)\cdots \phi_{l_{k(1)-1}}(z_{k(1)-1})
\\
\qquad{}\times[\dots[\phi_{l_{k(1)}}(z_{k(1)}), J^{-}_{\mu_{1, 1}}(u_{1, 1})]_{q^{l_{k(1)}}},
J^{-}_{\mu_{1, 2}}(u_{1, 2})]_{q^{l_{k(1)}-2}}\dots ,
J^{-}_{\mu_{1, n_1}}(u_{1, n_1})]_{q^{l_{k(1)}-2(n_1-1)}}
\dots
\\
\qquad{}\times[\dots[\phi_{l_{k(r)}}(z_{k(r)}), J^{-}_{\mu_{r, 1}}(u_{r, 1})]_{q^{l_{k(r)}}},
J^{-}_{\mu_{r, 2}}(u_{r, 2})]_{q^{l_{k(r)}-2}}\dots ,
J^{-}_{\mu_{r, n_r}}(u_{r, n_r})]_{q^{l_{k(r)}-2(n_r-1)}}
\\
\qquad{}\times
\phi_{l_{k(r)+1}}(z_{k(r)+1})\dots \phi_{l_n}(z_n)
S_{\epsilon_N}(t_N)\dots S_{\epsilon_1}(t_1)
|L\Big\rangle.
\end{gather*}
and the integrand in the right hand side signif\/ies to take the coef\/f\/icient of $\Bigg(\prod\limits_{\substack{1\le i \le r \\ 1\le j \le n_i}}u_{i, j}\Bigg)^{-1}$. For the notation $[x, y]_q$ see Appendix~\ref{appendixC}.

Let $(m)=(m_1, \dots, m_r)$, $0\le m_i \le n_i$.
%$\oint \left(\prod_{\genfrac{}{}{0pt}{}{1\le i_1 \le r}{1\le i_2 \le n_{i_1}}}\mu_{i_1, i_2}\frac{d u_{i_1, i_2}}{2\pi %iu_{i_1, i_2}}\right)
%F_{(\epsilon)(\mu)}^{(\nu)}(t, z)$ can be written as follows
Then
\begin{gather*}
\oint
\left(\prod_{\substack{1\le i_1 \le r \\ 1\le i_2 \le n_{i_1}}}\mu_{i_1, i_2}\frac{d u_{i_1, i_2}}{2\pi iu_{i_1, i_2}}\right)
F_{(\epsilon)(\mu)}^{(\nu)}(t, z) \\
\qquad{}=
\sum_{\genfrac{}{}{0pt}{}{0\le m_i\le n_i}{1\le i \le r}}
(-1)^{\sum\limits_{i=1}^{r} m_i}
\left(\prod_{i=1}^{r}q^{m_i l_{k(i)}}q^{-m_i(n_i-1)}\qbi{n_i}{m_i}\right)
\\
\qquad{}\times
\int_{C^N} \left(\prod_{\substack{1\le i_1 \le r \\ 1\le i_2\le n_{i_1}}}\mu_{i_1, i_2}\frac{d u_{i_1, i_2}}{2 \pi i u_{i_1, i_2}}\right)
F_{(\epsilon)(\mu)(m)}^{(\nu)}(t, z| u),
\end{gather*}
where
\begin{gather*}
F_{(\epsilon)(\mu)(m)}^{(\nu)}(t, z| u)
=\Big\langle L+\sum_{i=1}^{n}l_i-2N|\phi_{l_1}(z_1) \cdots
\phi_{l_{k(1)-1}}(z_{k(1)-1})
\nonumber
\\
\qquad{}\times
\big(J^{-}_{\mu_{1, 1}}(u_{1, 1}) \cdots J^{-}_{\mu_{1, m_1}}(u_{1, m_1}) \phi_{l_{k(1)}}(z_{k(1)})
J^{-}_{\mu_{1, m_1+1}}(u_{1, m_1+1}) \cdots J^{-}_{\mu_{1, n_1}}(u_{1, n_1})\big)\cdots
\\
\qquad{}\times
\big(J^{-}_{\mu_{r, 1}}(u_{r, 1}) \cdots J^{-}_{\mu_{r, m_r}}(u_{r, m_r}) \phi_{l_{k(r)}}(z_{k(r)})
J^{-}_{\mu_{r, m_r+1}}(u_{r, m_r+1}) \cdots J^{-}_{\mu_{r, n_r}}(u_{r, n_r})\big)
\\
\qquad{}\times
\phi_{l_{k(r)+1}}(z_{k(r)+1}) \cdots \phi_{l_{n}}(z_{n}) S_{\epsilon_N}(t_N) \cdots S_{\epsilon_1}(t_1)|L\Big\rangle,
\end{gather*}
and $C^N$ is a suitable deformation of the torus ${\mathbb{T}}^N$ specif\/ied as follows.
We introduce the lexicographical order
\begin{gather*}
(i_1, i_2)<(j_1, j_2) \ \ \Leftrightarrow   \ \  i_1<j_1 \quad \text{or} \quad i_1=j_1 \quad \text{and} \quad i_2<j_2.
\end{gather*}
For a given $(m)=(m_1, \dots, m_r)$, $1\le m_i \le n_i$, we def\/ine
\begin{gather*}
j<(i_1, i_2) \ \ \Leftrightarrow \ \ j<k(i_1) \quad \text{or} \quad j=k(i_1) \quad \text{and} \quad m_{i_1}<i_2,
\\
j>(i_1, i_2) \ \ \Leftrightarrow \ \ j>k(i_1) \quad \text{or} \quad j=k(i_1) \quad \text{and}\quad  m_{i_1}\ge i_2.
\end{gather*}
The contour for the integration variable $u_{i_1, i_2}$ is a simple closed curve rounding the origin in the counterclockwise direction such that $q^{l_j+k+2}z_j$ $((i_1, i_2)<j)$, $q^{-2}u_{j_1, j_2}$ $((i_1, i_2)<(j_1, j_2))$, $q^{-\mu_{i_1, i_2}(k+2)}t_a$ $(1\le a \le N)$ are inside, and
$q^{-l_j+k+2}z_j$ $((i_1, i_2)>j)$, $q^2 u_{j_1, j_2}$ $((j_1, j_2)<(i_1, i_2))$ are outside.
We denote it $C_{(i_1, i_2)}$.

Then
\begin{gather*}
F_{(\epsilon)(\mu)(m)}^{(\nu)}(t, z| u)=f^{(\nu)}(t, z)\Phi(t, z)G_{(\epsilon)(\mu)(m)}^{(\nu)}(t, z| u),
\end{gather*}
where
\begin{gather*}
f^{(\nu)}(t, z)
 = \left\{\prod_{i<j}(q^k z_i)^{\frac{l_i l_j}{2(k+2)}} \xi_{l_i, l_j}(z_i/z_j)\right\}
\left\{\prod_{i=1}^{n}(q^k z_i)^{-\frac{N l_i}{k+2}}\right\}
\\
\phantom{f^{(\nu)}(t, z)=}{} \times \left\{\prod_{i=1}^{n}(q^k z_i)^{\frac{L l_i}{2(k+2)}}\right\}  \left\{\prod_{i=1}^{N}(q^{-2}t_i)^{-\frac{L}{k+2}}\right\}
\left\{\prod_{a<b} (q^{-2} t_b)^{\frac{2}{k+2}}\right\},
\\
G_{(\epsilon)(\mu)(m)}^{(\nu)}(t, z| u)
 =
\widehat{G}_{(\epsilon)(\mu)(m)}^{(\nu)}(t, z| u)
\left(\prod_{a<b} \frac{q^{\epsilon_b} t_b-q^{\epsilon_a} t_a}{t_b-q^{-2} t_a}\right),
\\
\widehat{G}_{(\epsilon)(\mu)(m)}^{(\nu)}(t, z| u) = \left(\prod_{(i_1, i_2)}q^{L\mu_{i_1, i_2}}\right)
\left(\prod_{(i_1, i_2)>j}\frac{z_j-q^{\mu_{i_1, i_2} l_j-k-2}u_{i_1, i_2}}{z_j-q^{l_j-k-2}u_{i_1, i_2}}\right)
\\
\phantom{\widehat{G}_{(\epsilon)(\mu)(m)}^{(\nu)}(t, z| u) =}{}  \times
\left(\prod_{(i_1, i_2)<j}q^{\mu_{i_1, i_2} l_j}
\frac{u_{i_1, i_2}-q^{-\mu_{i_1, i_2} l_j+k+2}z_j}{u_{i_1, i_2}-q^{l_j+k+2}z_j}\right)
\\
\phantom{\widehat{G}_{(\epsilon)(\mu)(m)}^{(\nu)}(t, z| u) =}{} \times
\left(\prod_{\genfrac{}{}{0pt}{}{(i_1, i_2)}{1\le b \le N}}q^{-\mu_{i_1, i_2}}
\frac{u_{i_1, i_2}-q^{-\mu_{i_1, i_2}(k+1)-\epsilon_b} t_b}{u_{i_1, i_2}-q^{-\mu_{i_1, i_2}(k+2)} t_b}\right)
 \\
\phantom{\widehat{G}_{(\epsilon)(\mu)(m)}^{(\nu)}(t, z| u) =}{} \times
\left(\prod_{(i_1, i_2)<(j_1, j_2)}
\frac{q^{-\mu_{i_1, i_2}} u_{i_1, i_2}-q^{-\mu_{j_1, j_2}}u_{j_1, j_2}}
{u_{i_1, i_2}-q^{-2}u_{j_1, j_2}}\right).
\end{gather*}

For $i$, let $A_{\mu, i}^{\pm}=\{(i, j) |\mu_{i, j}=\pm \}$. The number of elements in $A_{\mu, i}^{\pm}$
is $a_i^{\pm}$and $A_{\mu, i}^{\pm}=\{{{\ell}}^{\pm}_{i, 1}, \dots, {{\ell}}^{\pm}_{i, a_i^{\pm}}\}$ .
We set $a_i^{-}=a_i$, $A_{\mu, i}^{-}=A_{\mu, i}$, $A_{\mu}=\cup_{i=1}^{r}A_{\mu, i}$ and
\begin{gather*}
\widehat{J}_{(\epsilon)(\mu)}^{(\nu)} = \sum_{\substack{0\le m_i\le n_i \\ 1\le i \le r}}(-1)^{\sum\limits_{i_1}^{r}m_i}
\left\{\prod_{i=1}^{r}q^{m_i l_{k(i)}}q^{-m_i(n_i-1)}\qbi{n_i}{m_i}\right\}
\\
\phantom{\widehat{J}_{(\epsilon)(\mu)}^{(\nu)} =}{}
\times\int_{C^N}\left(\prod_{(i_1,i_2)}\mu_{i_1, i_2}\frac{d u_{i_1, i_2}}{2\pi i u_{i_1, i_2}}\right)
\widehat{G}_{(\epsilon) (\mu)(m)}^{(\nu)}.
\end{gather*}
See the beginning of the next section for the notation of the $q$-binomial coefficient $\qbi{n_i}{m_i}$.

For a given $(a)=(a_1, \dots, a_r)$, $1\le a_i \le n_i$, we def\/ine $\widehat{J}_{(\epsilon)(a)}^{(\nu)}$ and ${J}_{(a)}^{(\nu)}$ as follows
\begin{gather*}
 \widehat{J}_{(\epsilon)(a)}^{(\nu)}
=\sum_{\substack{|A_{\mu, i}|=a_i \\ 1\le i \le r}}\widehat{J}_{(\epsilon)(\mu)}^{(\nu)},
\\
 {J}_{(a)}^{(\nu)}=\sum_{\epsilon_1, \dots, \epsilon_N=\pm}\left(\prod_{j=1}^N \epsilon_j\right)
\left(\prod_{1\le a <b \le N}\frac{q^{\epsilon_b}t_b-q^{\epsilon_a}t_a}{t_b-q^{-2}t_a}\right)
\widehat{J}_{(\epsilon)(a)}^{(\nu)}.
\end{gather*}

Using $J^{(\nu)}_{(a)}$, $F^{(\nu)}(t, z)$ can be written as
\begin{gather*}
F^{(\nu)}(t, z)=(-1)^N\big(q-q^{-1}\big)^{-2N}\left(\prod_{i=1}^r\frac{1}{[n_i]!}\right)\left(\prod_{b=1}^N t_b^{-1}\right)
f^{(\nu)}(t, z)\Phi(t,z)\sum_{(a)}J^{(\nu)}_{(a)}.
\end{gather*}
Theorem~\ref{main} straightforwardly follows from the following proposition.

\begin{proposition}\label{mainprop}
If $(a)\ne(n_1, n_2, \dots, n_r)$, $J_{(a)}^{(\nu)}(t, z)=0$. For $(a)=(n_1, n_2, \dots, n_r)$ we have
\begin{gather*}
J_{(n_1, \dots, n_r)}^{(\nu)}(t, z) =
(-1)^N\big(1-q^{-2}\big)^{N} q^{N(N-L)+\frac{N(N-1)}{2}-\big(\sum\limits_{i=1}^n l_i\big)N}
\\
\phantom{J_{(n_1, \dots, n_r)}^{(\nu)}(t, z) = }{}
\times
\prod_{s=1}^{r}\left\{q^{\big(\sum\limits_{t=s+1}^r n_t\big) n_s -l_{k(s)}n_s}[n_s]!
\prod_{i=0}^{n_s-1}(1-q^{2(l_{k(s)}-i)})\right\}w_{(-\nu)}(t, z).
\end{gather*}
\end{proposition}
This proposition is proved by performing integrals in the variables $u_{i, j}$ in the next section.

\section{Proof of Proposition \ref{mainprop}}\label{section5}
We set
\begin{gather*}
[n]! = \prod_{i=1}^{n}[i],
\qquad
\qbi{n}{m} = \frac{[n]!}{[n-m]![m]!},
\end{gather*}
for nonnegative integers $n$, $m$ $(n\ge m)$.
To prove Proposition~\ref{mainprop}, we have to calculate $\widehat{J}^{(\nu)}_{(\epsilon)(a)}$.
We need the following lemmas.
\begin{lemma}{\label{combi}}
For $n\ge 1$ and $n\ge m \ge 0$, we have
\begin{gather*}
(i)~ \ \ \sum_{\substack{A\sqcup B=\{1, 2, \dots, n\}\\ |A|=m}}
\left(\prod_{\substack{i<j \\ i \in A , j \in B}} q^2 \right)
=
q^{m(n-m)}\qbi{n}{m};
\\
(ii) \ \ \sum_{\substack{A\sqcup B=\{1, 2, \dots, n\}\\  |A|=m \\ \mu_i=1 (i\in A), \, \mu_i=-1 (i\in B)}}  \left(\prod_{i<j} q^{\mu_i} \right)
=
q^{-\frac{n(n-1)}{2}+m(n-1)}\qbi{n}{m}.
\end{gather*}
\end{lemma}

\begin{proof} By the $q$-binomial theorem
\begin{gather*}
\prod_{i=1}^n\big(1+q^{-n-1+2i}x\big)=\sum_{i=0}^n\qbi{n}{i}x^i,
\end{gather*}
we have the equation
\begin{gather*}
\sum_{1\le i_1<\dots <i_m\le n}q^{2\sum\limits_{j=1}^{m}i_j}=q^{(n+1)m}\qbi{n}{m}.
\end{gather*}
The assertions $(i)$ and $(ii)$ easily follow from this equation.
\end{proof}

\begin{lemma}{\label{sym}}
Let $n\ge 1$, $n\ge m\ge 0$ and $1\le i_1<\dots<i_m\le n$. Then we have
\begin{gather*}
 \sum_{\sigma\in S_n}{\rm{sgn}}\,\sigma \,\,t_{\sigma(i_1)} t_{\sigma(i_2)} \cdots t_{\sigma(i_m)}
\prod_{1\le a<b\le n}(t_{\sigma(b)}-q^{-2}t_{\sigma(a)})
\\
\qquad{} =
q^{-m(n+1)-\frac{n(n-1)}{2}+2\sum\limits_{j=1}^{m}i_j}[m]! [n-m]!\,\,e_m(t_1, \dots, t_n)
\prod_{1\le a<b\le n}(t_b-t_a) ,
\end{gather*}
where $e_m(t_1, \dots, t_n)$ is the $m$-th elementary symmetric polynomial.
\end{lemma}

\begin{proof} Set
\begin{gather*}
F(t)=\sum_{\sigma\in S_n}{\rm{sgn}}\,\sigma \,\,t_{\sigma(i_1)} t_{\sigma(i_2)} \cdots t_{\sigma(i_m)}
\prod_{1\le a<b\le n}(t_{\sigma(b)}-q^{-2}t_{\sigma(a)}).
\end{gather*}

It is easy to see that $F(t)$ is an antisymmetric polynomial. So we can write
\[
F(t)=S(t)\prod_{1\le a<b\le n}(t_b-t_a),
 \]
 where $S(t)$  is a symmetric polynomial. Moreover $S(t)$ is a homogeneous polynomial of degree~$m$ and ${\rm{deg}}_{t_i} S(t)=1$ for all $i\in\{1, \dots , n\}$. Hence we have
 \[
 S(t)=c e_m(t)
 \]
  for some constant $c$.

The number $(-1)^{\sum\limits_{j=1}^{m}i_j+\frac{n(n-1)}{2}-\frac{m(m+1)}{2}}c$ is equal to the coef\/f\/icient of
\[
t_{i_1}^{n}t_{i_2}^{n-1}\cdots
 t_{i_m}^{n-m+1}t_1^{n-m-1}t_{2}^{n-m-2}\cdots t_{n-1}
\] in $F(t)$.

We can show
\begin{gather*}
c=q^{-2nm+m(m-1)+2\sum\limits_{k=1}^{m} i_k}
\left(q^{-m(m-1)}\sum_{\sigma\in S_{m}}q^{2\ell(\sigma)}\right)
\left(q^{-(n-m)(n-m-1)}\sum_{\tau\in S_{n-m}}q^{2\ell(\tau)}\right),
\end{gather*}
where $\ell(\sigma)$ is  the inversion number of $\sigma$.

Using the fact $\sum\limits_{\sigma\in S_{m}}q^{2\ell(\sigma)}=q^{\frac{m(m-1)}{2}}[m]!$, we have the desired result.
\end{proof}

\begin{lemma}{\label{z}}
For $1\le n \le l$, we have
\begin{gather*}
\sum_{s=0}^n(-1)^s q^{-s(n-1)}
\qbi{n}{s}
\sum_{\sigma \in S_n}
\prod_{i=1}^{s}(z-q^{l}t_{\sigma(i)})
\prod_{i=s+1}^{n}\big(z-q^{-l}t_{\sigma(i)}\big)
\left(\prod_{1\le a <b \le n}\!\!\frac{t_{\sigma(b)}-q^{-2}t_{\sigma(a)}}{t_{\sigma(b)}-t_{\sigma(a)}}\right)
\\
\qquad{} =
(-1)^n q^{-ln-\frac{n(n-1)}{2}}
\left\{\prod_{i=0}^{n-1}\big(1-q^{2(l-i)}\big)\right\} [n]! t_1t_2\dots t_n.
\end{gather*}
\end{lemma}

\begin{proof}
We set
\begin{gather*}
L_{n, s} =
\sum_{\sigma \in S_n} {\rm{sgn}}\,\sigma  \prod_{i=1}^{s}(z-q^{l}t_{\sigma(i)})
\prod_{j=s+1}^{n}\big(z-q^{-l}t_{\sigma(j)}\big)
\left(\prod_{i>j}\frac{t_{\sigma(i)}-q^{-2}t_{\sigma(j)}}{t_{i}-t_{j}}\right),
\\
L_n=
\sum_{s=0}^n(-1)^s q^{-s(n-1)} \qbi{n}{s} L_{n, s}.
\end{gather*}
Using Lemma \ref{sym},
\begin{gather*}
L_{n, s} = \sum_{k=0}^n (-1)^k z^{n-k} e_k(t) q^{-k(n+1)-\frac{n(n-1)}{2}}
[k]![n-k]!
\!\left\{\sum_{t=0}^k q^{2lt-lk}
\!\left(\!\sum_{\substack{1\le i_1< i_2<\cdots< i_t\le s \\ s< i_{t+1}<\cdots<i_{k}\le n}}\!\!q^{\sum\limits_{j=1}^k 2 i_j}\right)\!
\right\}
\\
\phantom{L_{n, s}}{} =\sum_{k=0}^n (-1)^k z^{n-k} e_k(t) q^{-lk-k(n+1)-\frac{n(n-1)}{2}} [k]![n-k]!
\\
\phantom{L_{n, s}=}{}\times
\left(\sum_{t=0}^k q^{2lt}
q^{2s(k-t)+(s+1)t+(n-s+1)(k-t)}\qbi{s}{t}\qbi{n-s}{k-t}
\right).
\end{gather*}
Then,
\begin{gather*}
L_n = \sum_{s=0}^n (-1)^s q^{-s(n-1)} \qbi{n}{s}
\sum_{k=0}^n (-1)^k z^{n-k} e_k(t)q^{-lk-k(n+1)-\frac{n(n-1)}{2}}  [k]![n-k]!
\\
\phantom{L_n =}{} \times
\left(\sum_{t=0}^k q^{2lt}q^{sk+k+n(k-t)}
\qbi{s}{t}\qbi{n-s}{k-t}
\right)
\\
\phantom{L_n}{}=
[n]!\sum_{k=0}^n (-1)^k z^{n-k} e_k(t)\, q^{-lk-k(n+1)-\frac{n(n-1)}{2}}
\\
\phantom{L_n =}{}\times
\sum_{t=0}^k q^{2lt}q^{(k-t)(n+1)+t}  \qbi{k}{t}
\sum_{s=t}^{n-k+t} (-1)^s q^{-s(n-k-1)} \qbi{n-k}{s-t}
\\
\phantom{L_n}{}=[n]!\sum_{k=0}^n (-1)^k z^{n-k} e_k(t)\, q^{-lk-k(n+1)-\frac{n(n-1)}{2}}\\
\phantom{L_n =}{}\times
\sum_{t=0}^k q^{2lt} q^{(k-t)(n+1)+t} \qbi{k}{t} (-1)^{t}q^{-t(n-k-1)}\delta_{n, k}
\\
\phantom{L_n}{}=
[n]!(-1)^n q^{-ln} q^{-\frac{n(n-1)}{2}}
\sum_{t=0}^n  (-1)^{t}q^{2lt} q^{-(n-1)t} \qbi{n}{t}\, e_n(t)
\\
\phantom{L_n}{}=[n]!(-1)^n q^{-ln} q^{-\frac{n(n-1)}{2}}
\left\{\prod_{i=0}^{n-1}(1-q^{2(l-i)})\right\}  e_n(t).
\end{gather*}
Here we have used the $q$-binomial theorem.
\end{proof}

 For a given sequence $(m_i)_{i=1}^{r}$ $(0\le m_i \le n_i)$, let $M_i=\{(i, j)\,|\,j\le m_i\}$.
Set
\begin{gather*}
\widehat{I}_{(\mu)(\epsilon)(m)}^{(\nu)}
=\int_{C^N}\left(\prod_{(i_1,i_2)}\frac{d u_{i_1, i_2}}{2\pi i u_{i_1, i_2}}\right)
 \widehat{G}_{(\mu)(\epsilon)(m)}^{(\nu)}.
\end{gather*}

\begin{lemma} \label{i}
We have
\begin{gather*}
\widehat{I}_{(\mu)(\epsilon)(m)}^{(\nu)}
=q^{(L-N)\big\{\sum\limits_{s=1}^r(n_s-2a_s)\big\} }
\left(\prod_{(i_1, i_2)<j} q^{\mu_{i_1, i_2} l_j}\right)
\left(\prod_{(i_1, i_2)<(j_1, j_2)} q^{-\mu_{i_1, i_2}}\right)
\\
\times
 \sum_{\substack{C_i\sqcup D_i=A_{\mu, i}\\ D_i'=D_i\cap M_i \\ 1\le i \le r}}
\left(\prod_{b=1}^{N}q^{-1-\epsilon_b}\right)^{\sum\limits_{i=1}^{r}|C_i|}\left(\prod_{\substack{(i_1, i_2)<(j_1, j_2)\\ (i_1, i_2)\in C_1\cup\cdots \cup C_r \\ (j_1, j_2)\in D_1\cup\cdots \cup D_r}} q^2\right)
\\
\times
\sum_{\substack{1\le b_{i, j} \le N \\ 1\le i \le r \\ 1\le j \le |D_i|}}
\prod_{i_1=1}^r \left\{
\prod_{i_2=1}^{|D_{i_1}|}\left((1-q^{-1-\epsilon_{b_{i_1, i_2}}})
\prod_{b\ne b_{i_1, i_2}}\frac{t_{b_{i_1, i_2}}-q^{-1-\epsilon_b}t_b}{t_{b_{i_1, i_2}}-t_b}
\right.
\right.
\\
\times
\left.
\left.
\prod_{(i_1, i_2)<(j_1, j_2)}\frac{t_{b_{i_1, i_2}}-t_{b_{j_1, j_2}}}{t_{b_{i_1, i_2}}-q^{-2}t_{b_{j_1, j_2}}}
\prod_{j=1}^{k(i_1)-1}\frac{z_j-q^{-l_j}t_{b_{i_1, i_2}}}{z_j-q^{l_j}t_{b_{i_1, i_2}}}\right)
\prod_{i_2=|D_{i_1}'|+1}^{|D_{i_1}|}\frac{z_{k(i_1)}-q^{-l_{k(i_1)}}t_{b_{i_1, i_2}}}{z_{k(i_1)}-q^{l_{k(i_1)}}t_{b_{i_1, i_2}}}
\right\}.
\end{gather*}
\end{lemma}

\begin{proof}
We integrate with respect the variables $u_{i, j}$, $(i, j)\in A_{\mu}^+$ in the order
$u_{{{\ell}}^+_{1, 1}}, \dots, u_{{{\ell}}^+_{1, a_1^+}},$
$ u_{{{\ell}}^+_{2, 1}}, \dots, u_{{{\ell}}^+_{r, a_r^+}}$.
With respect to $u_{{\bf{\ell}}_{1, 1}^ +}$
the only singularity outside $C_{{\bf{\ell}}^+_{1, 1}}$ is~$\infty$.
Then the integral in $u_{{\bf{\ell}}_{1, 1}^ +}$ is calculated by taking the residue at $\infty$. After this integration the integrand as a function of $u_{{{\ell}}^+_{1, 2}}$ has a similar
structure. Then the integral with respect to $u_{{{\ell}}^+_{1, 2}}$ is calculated by taking residue at $\infty$ and so~on. Finally we get
\begin{gather*}
\widehat{I}^{(\nu)}_{(\epsilon)(\mu)(m)} = (-1)^{\sum\limits_{i=1}^{r}a_i^{+}}
{{\mathop{\rm Res}_{{u_{{{\ell}}_{r, a_r^{+}}^{+}}=\infty}}}}\cdots{{\mathop{\rm Res}_{{u_{{{\ell}}_{r, 1}^{+}}=\infty}}}}\cdots
{{\mathop{\rm Res}_{{u_{{{\ell}}_{1, a_1^{+}}^{+}}=\infty}}}}\cdots{{\mathop{\rm Res}_{{u_{{{\ell}}_{1, 1}^{+}}=\infty}}}}
\widehat{G}_{(\epsilon)(\mu)(m)}^{\nu}(t, z| u)
\\
\phantom{\widehat{I}^{(\nu)}_{(\epsilon)(\mu)(m)} }{} =
\left(\prod_{(i_1, i_2)}q^{(L-N)\mu_{i_1, i_2}}\right)
\left(\prod_{(i_1, i_2)<j}q^{\mu_{i_1, i_2} l_j}\right)
\left(\prod_{(i_1, i_2)<(j_1, j_2)}q^{-\mu_{i_1, i_2}}\right)
\\
\phantom{\widehat{I}^{(\nu)}_{(\epsilon)(\mu)(m)} =}{}\times
\int_{C^{N-\sum\limits_{i=1}^{r}a_i^+}}
\left(\prod_{(i_1, i_2)\in A_{\mu}}\frac{d u_{i_1, i_2}}{2 \pi i u_{i_1, i_2}}\right)
\left(\prod_{\substack{j<(i_1, i_2)\\ (i_1, i_2)\in A_{\mu}}}\frac{z_j-q^{-l_j-k-2}u_{i_1, i_2}}{z_j-q^{l_j-k-2}u_{i_1, i_2}}\right)
\\
\phantom{\widehat{I}^{(\nu)}_{(\epsilon)(\mu)(m)} =}{}\times
\left(\prod_{\genfrac{}{}{0pt}{}{(i_1, i_2)\in A_{\mu}}{1\le b \le N}}
\frac{u_{i_1, i_2}-q^{k+1-\epsilon_b} t_b}{u_{i_1, i_2}-q^{k+2} t_b}\right)
\left(\prod_{\substack{(i_1, i_2)<(j_1, j_2)\\ (i_1, i_2), (j_1, j_2)\in A_{\mu}}}
\frac{u_{i_1, i_2}-u_{j_1, j_2}}
{u_{i_1, i_2}-q^{-2}u_{j_1, j_2}}\right),
\end{gather*}
where $C^{N-\sum\limits_{i=1}^{r}a_i^+}$
is the resulting contour for $(u_{{\bf{\ell}}_{1, 1}}, \dots, u_{{\bf{\ell}}_{r, a_r}})$.
We set
\begin{gather*}
I_{(\epsilon)(\mu)(m)}^{(\nu) +}(t, z) = \left(\prod_{(i_1, i_2)\in A_{\mu}}\frac{1}{u_{i_1, i_2}}\right)
\left(\prod_{\substack{j<(i_1, i_2)\\ (i_1, i_2)\in A_{\mu}}}
\frac{z_j-q^{-l_j-k-2}u_{i_1, i_2}}{z_j-q^{l_j-k-2}u_{i_1, i_2}}\right)
\\
\phantom{I_{(\epsilon)(\mu)(m)}^{(\nu) +}(t, z)}{}
\times
\left(\prod_{\substack{(i_1, i_2)\in A_{\mu} \\ 1\le b \le N}}
\frac{u_{i_1, i_2}-q^{k+1-\epsilon_b} t_b}{u_{i_1, i_2}-q^{k+2} t_b}\right)
\left(\prod_{\substack{(i_1, i_2)<(j_1, j_2) \\ (i_1, i_2), (j_1, j_2)\in A_{\mu}}}
\frac{u_{i_1, i_2}-u_{j_1, j_2}}
{u_{i_1, i_2}-q^{-2}u_{j_1, j_2}}\right).
\end{gather*}
Next we perform integrations with respect to the remaining variables $u_{i, j}$, $(i, j)\in A_{\mu}$ in the order
$u_{{{\ell}}_{r, a_r}}, \dots$, $u_{{{\ell}}_{r, 1}}, u_{{{\ell}}_{r-1, a_{r-1}}}, \dots, u_{{{\ell}}_{1, 1}}$.
The poles of the integrand inside $C_{{{\ell}}_{r, a_r}}$ are $0$ and $q^{k+2}t_b$, $b=1,\dots, N$. Thus we have
\begin{gather*}
\int_{C_{{{\ell}}_{r, a_r}}}\frac{d u_{{{\ell}}_{r, a_r}}}{2 \pi i } I_{(\epsilon)(\mu)(m)}^{(\nu) +}(t, z)
\\
=
\left(\prod_{i\le b \le N}q^{-1-\epsilon_b}\right)
\left(\prod_{\substack{(i_1, i_2)\in A_{\mu}\\ (i_1, i_2)\ne {{\ell}}_{r, a_r}}}\frac{1}{u_{i_1, i_2}}\right)
\left(\prod_{\substack{j<(i_1, i_2)\\ (i_1, i_2)\in A_{\mu}-\{{{\ell}}_{r, a_r}\}}}\frac{z_j-q^{-l_j-k-2}u_{i_1, i_2}}{z_j-q^{l_j-k-2}u_{i_1, i_2}}\right)
\\
\times
\left(\prod_{\substack{(i_1, i_2)\in A_{\mu}-\{{{\ell}}_{r, a_r}\}\\ 1\le b \le N}}
\frac{u_{i_1, i_2}-q^{k+1-\epsilon_b} t_b}{u_{i_1, i_2}-q^{k+2} t_b}\right)
\left(\prod_{\substack{(i_1, i_2)<(j_1, j_2)<{{\ell}}_{r, a_r}\\ (i_1, i_2), (j_1, j_2)\in A_{\mu}}}
\frac{u_{i_1, i_2}-u_{j_1, j_2}}
{u_{i_1, i_2}-q^{-2}u_{j_1, j_2}}\right)
\\
+
\sum_{1\le b_{{{\ell}}_{r, a_r}}\le N}
(1-q^{-1-\epsilon_{b_{{{\ell}}_{r, a_r}}}})
\left(\prod_{j<{{\ell}}_{r, a_r}}\frac{z_j-q^{-l_j}t_{b_{{{\ell}}_{r, a_r}}}}{z_j-q^{l_j-k-2}t_{b_{{{\ell}}_{r, a_r}}}}\right)
\left(\prod_{\substack{1\le b \le N\\ b\ne b_{{{\ell}}_{r, a_r}}}}
\frac{t_{b_{{{\ell}}_{r, a_r}}}-q^{-1-\epsilon_b} t_b}{t_{b_{{{\ell}}_{r, a_r}}}-t_b}\right)
\\
\times
\left(\prod_{(i_1, i_2)<{{\ell}}_{r, a_r}}\!\!
\frac{u_{i_1, i_2}-q^{k+2}t_{b_{{{\ell}}_{r, a_r}}}}{u_{i_1, i_2}-q^{k}t_{b_{{{\ell}}_{r, a_r}}}}\right)
\left(\prod_{\substack{(i_1, i_2)\in A_{\mu}\\ (i_1, i_2)\ne {{\ell}}_{r, a_r}}}\frac{1}{u_{i_1, i_2}}\right)\!
\left(\! \prod_{\substack{j<(i_1, i_2) \\ (i_1, i_2)\in A_{\mu}-\{{{\ell}}_{r, a_r}\}}}\!\!\!\!\! \frac{z_j-q^{-l_j-k-2}u_{i_1, i_2}}{z_j-q^{l_j-k-2}u_{i_1, i_2}}\!\right)\!
\\
\times
\left(\prod_{\substack{(i_1, i_2)\in A_{\mu}-\{{{\ell}}_{r, a_r}\}\\ 1\le b \le N}}
\frac{u_{i_1, i_2}-q^{k+1-\epsilon_b} t_b}{u_{i_1, i_2}-q^{k+2} t_b}\right)
\left(\prod_{\substack{(i_1, i_2)<(j_1, j_2)<{{\ell}}_{r, a_r} \\ (i_1, i_2), (j_1, j_2)\in A_{\mu}}}
\frac{u_{i_1, i_2}-u_{j_1, j_2}}
{u_{i_1, i_2}-q^{-2}u_{j_1, j_2}}\right).
\end{gather*}
The integrand in $u_{{{\ell}}_{r, a_r-1}}$ has the poles at $0$ and $q^{k+2}t_b$ inside $C_{{{\ell}}_{r, a_r-1}}$ and so on.
Finally we get
\begin{gather*}
\widehat{I}^{(\nu)}_{(\epsilon)(\mu)}=\left(\prod_{(i_1, i_2)} q^{(L-N)\mu_{i_1, i_2} }\right)
\left(\prod_{(i_1, i_2)<j} q^{\mu_{i_1, i_2} l_j}\right)
\left(\prod_{(i_1, i_2)<(j_1, j_2)} q^{-\mu_{i_1, i_2}}\right)
\\
\phantom{\widehat{I}^{(\nu)}_{(\epsilon)(\mu)}=}{}
\times
\sum_{\substack{w_{{{\ell}}_{i_1, i_2}}\in \{0\}\cup(T-W_{i_1, i_2})\\ (i_1, i_2)\in A_{\mu}}}
\mathop{\rm Res}_{u_{{{\ell}}_{1, 1}}=w_{{{\ell}}_{1, 1}}}\cdots
\mathop{\rm Res}_{u_{{{\ell}}_{r, a_r}}=w_{{{\ell}}_{r, a_r}}}
I_{(\epsilon)(\mu)}^{(\nu)+},
\end{gather*}
where $T=\{t_1, t_2, \dots, t_N\}$,
${W_{i_1, i_2}=\mathop{\cup}\limits_{{{\ell}}_{i_1, i_2}<{{\ell}}_{j_1, j_2}}\{w_{{{\ell}}_{j_1, j_2}}\}}$.

Set $C_i=\{{{\ell}}_{i, j}\,|\, w_{{{\ell}}_{i, j}}=0\}$, $D_i=A_{\mu, i}-C_i$. Then we have the desired result.
\end{proof}

Now we can calculate $\widehat{J}_{(\epsilon)(a)}^{(\nu)}$.
\begin{proposition}
We have
\begin{gather*}
\widehat{J}_{(\epsilon)(a)}^{(\nu)}=
(-1)^{\sum\limits_{i=1}^r a_i}
\left(q^{\sum\limits_{s=1}^r\big(\sum\limits_{t=k(s)+1}^{n}l_t\big)(n_s-2 a_s)} \right)
\left(q^{(L-N)\big\{\sum\limits_{s=1}^{r}(n_s-2 a_s)\big\}}\right)
\\
\phantom{\widehat{J}_{(\epsilon)(a)}^{(\nu)}=}{}\times
\left(q^{-\sum\limits_{s=1}^{r}\sum\limits_{t=s+1}^r n_s(n_t-2 a_t)}\right)
\sum_{\substack{1\le b_{i_1, i_2}\le N \\ 1\le i_1 \le r \\ 1\le i_2 \le a_{i_1}}}
\left(\prod_{i_1<j_1}\frac{t_{b_{i_1, i_2}}-t_{b_{i_1, i_2}}}{t_{b_{i_1, i_2}}-q^{-2}t_{b_{i_1, i_2}}}\right)
\\
\phantom{\widehat{J}_{(\epsilon)(a)}^{(\nu)}=}{}\times
\prod_{i_1=1}^{r}\left\{\sum_{s_{i_1}=0}^{a_{i_1}}
q^{a_{i_1}(n_{i_1}-s_{i_1}-1)-\frac{n_{i_1}(n_{i_1}-1)}{2}}
\frac{[n_{i_1}]!}{[s_{i_1}]![a_{i_1}-s_{i_1}]!}
\right.
\\
\left.\phantom{\widehat{J}_{(\epsilon)(a)}^{(\nu)}=}{}\times
\sum_{i=0}^{n_{i_1}-a_{i_1}}
(-1)^{i+s_{i_1}}q^{i(2l_{k(i_1)}-n_{i_1}-a_{i_1}+1)+s_{i_1}}
\frac{1}{[i]![n_{i_1}-a_{i_1}-i]!}
\right.
\\
\left.\phantom{\widehat{J}_{(\epsilon)(a)}^{(\nu)}=}{}\times
\left\{\prod_{i_2=1}^{a_{i_1}}\left(\big(1-q^{-1-\epsilon_{b_{i_1, i_2}}}\big)
\prod_{b\ne b_{i_1, i_2}}\frac{t_{b_{i_1, i_2}}-q^{-1-\epsilon_b}t_b}{t_{b_{i_1, i_2}}-t_b}
\prod_{i_2<j_2}\frac{t_{b_{i_1, i_2}}-t_{b_{i_1, j_2}}}{t_{b_{i_1, i_2}}-q^{-2}t_{b_{i_1, j_2}}}
\right.
\right.
\right.
\\
\left.
\left.
\left.
\phantom{\widehat{J}_{(\epsilon)(a)}^{(\nu)}=}{}
\times
\prod_{j=1}^{k(i_1)-1}\frac{z_j-q^{-l_j}t_{b_{i_1, i_2}}}{z_j-q^{l_j}t_{b_{i_1, i_2}}}
\right)
\prod_{i_2=s_{i_1}+1}^{a_{i_1}}\frac{z_{k(i_1)}-q^{-l_{k(i_1)}}t_{b_{i_1, i_2}}}{z_{k(i_1)}-q^{l_{k(i_1)}}t_{b_{i_1, i_2}}}
\right\}
\right\}.
\end{gather*}
\end{proposition}

\begin{proof} Using Lemma~\ref{i} we have
\begin{gather}
%1
\widehat{J}_{(\epsilon)(a)}^{(\nu)} = (-1)^{\sum\limits_{i=1}^r a_i}\sum_{\substack{|A_{\mu, i}|=a_i \\ 1\le i \le r}}
\sum_{\substack{0\le m_i\le n_i \\ 1\le i \le r}}(-1)^{\sum\limits_{i=1}^{r}m_i}
\left\{\prod_{i=1}^{r}q^{m_i l_{k(i)}}q^{-m_i(n_i-1)}\qbi{n_i}{m_i}\right\}
\nonumber
\\
\phantom{\widehat{J}_{(\epsilon)(a)}^{(\nu)} =}{}
\times
\left(\prod_{i=1}^rq^{(L-N)(n_i-2a_i)}\right)
\left(\prod_{(i_1, i_2)<j} q^{\mu_{i_1, i_2} l_j}\right)
\left(\prod_{(i_1, i_2)<(j_1, j_2)} q^{-\mu_{i_1, i_2}}\right)
\nonumber
\\
\phantom{\widehat{J}_{(\epsilon)(a)}^{(\nu)} =}{}
\times
\sum_{\substack{C_i\sqcup D_i=A_{\mu, i} \\ D_i'=D_i\cap M_i \\ 1\le i \le r}}
(\prod_{b=1}^{N}q^{-1-\epsilon_b})^{\sum\limits_{i=1}^{r}|C_i|}\left(\prod_{\substack{(i_1, i_2)<(j_1, j_2)\\ (i_1, i_2)\in C_1\cup\dots \cup C_r \\ (j_1, j_2)\in D_1\cup\cdots \cup D_r}} q^2\right)
\nonumber
\\
\phantom{\widehat{J}_{(\epsilon)(a)}^{(\nu)} =}{}
\times
\sum_{\substack{1\le b_{i, j} \le N \\ 1\le i \le r \\ 1\le j \le |D_i|}}
\prod_{i_1=1}^r \left\{
\prod_{i_2=1}^{|D_{i_1}|}\left(\big(1-q^{-1-\epsilon_{b_{i_1, i_2}}}\big)
\prod_{b\ne b_{i_1, i_2}}\frac{t_{b_{i_1, i_2}}-q^{-1-\epsilon_b}t_b}{t_{b_{i_1, i_2}}-t_b}
\right.
\right.
\nonumber
\\
\left.
\phantom{\widehat{J}_{(\epsilon)(a)}^{(\nu)} =}{}
\times
\prod_{(i_1, i_2)<(j_1, j_2)}\frac{t_{b_{i_1, i_2}}-t_{b_{j_1, j_2}}}{t_{b_{i_1, i_2}}-q^{-2}t_{b_{j_1, j_2}}}
\prod_{j=1}^{k(i_1)-1}\frac{z_j-q^{-l_j}t_{b_{i_1, i_2}}}{z_j-q^{l_j}t_{b_{i_1, i_2}}}
\right)
\nonumber
\\
\left.
\phantom{\widehat{J}_{(\epsilon)(a)}^{(\nu)} =}{}\times
\prod_{i_2=|D_{i_1}'|+1}^{|D_{i_1}|}\frac{z_{k(i_1)}-q^{-l_{k(i_1)}}t_{b_{i_1, i_2}}}{z_{k(i_1)}-q^{l_{k(i_1)}}t_{b_{i_1, i_2}}}
\right\}.
\label{rep}
\end{gather}
Set $\lambda_i=|A_{\mu, i}\cap M_i|$, $\gamma_i=|D_i|$, $s_i=|D_i'|$, $1\le i \le r$.
Then the right hand side of (\ref{rep}) is equal to
\begin{gather*}
 \sum_{\substack{0\le m_i\le n_i \\ 1\le i \le r}}(-1)^{\sum\limits_{i=1}^{r} m_i} \left\{\prod_{i=1}^{r}q^{m_i l_{k(i)}}q^{-m_i(n_i-1)}
\qbi{n_i}{m_i}\right\}
\\
\qquad{}\times\sum_{\substack{0\le \gamma_j \le a_j \\ 1\le j \le r}}\sum_{\substack{0\le s_j\le\gamma_j \\ 1\le j \le r}}
\sum_{\substack{0\le \lambda_j\le m_j \\ 1\le j \le r}}
C_{(a)(\gamma)}
\left\{q^{\sum\limits_{s=1}^{r}l_{k(s)}(m_s-2\lambda_s)}\right\}
\\
\qquad{}\times
\left\{\prod_{i_1=1}^r \left(\sum_{\substack{|A_{\mu, i_1}|=a_{i_1}\\ |A_{\mu, i_1}\cap M_{i_1}|=\lambda_{i_1}}}\prod_{\substack{i_2<j_2 \\ i_1=j_1}}
q^{-\mu_{i_1, i_2}}\right)\right\}\\
\qquad{}\times
\left(\prod_{i_1=1}^{r}
q^{\lambda_{i_1}\gamma_{i_1}+a_{i_1}\gamma_{i_1}-a_{i_1}s_{i_1}-\gamma_{i_1}^2}
\qbi{\lambda_{i_1}}{s_{i_1}}
\qbi{a_{i_1}-\lambda_{i_1}}{\gamma_{i_1}-s_{i_1}}
\right)
\\
\qquad{}\times
\sum_{\substack{1\le b_{i_1, i_2}\le N \\ 1\le i_1 \le r \\ 1\le i_2 \le \gamma_{i_1}}}
\prod_{i_1=1}^{r}\left\{\left(\prod_{i_2=1}^{\gamma_{i_1}}(1-q^{-1-\epsilon_{b_{i_1, i_2}}})
\prod_{b\ne b_{i_1, i_2}}\frac{t_{b_{i_1, i_2}}-q^{-1-\epsilon_b}t_b}{t_{b_{i_1, i_2}}-t_b}
\right.
\right.
\\
\left.
\left.
\qquad{}\times
\prod_{i_2<j_2}\frac{t_{b_{i_1, i_2}}-t_{b_{i_1, j_2}}}{t_{b_{i_1, i_2}}-q^{-2}t_{b_{i_1, j_2}}}
\prod_{j=1}^{k(i_1)-1}\frac{z_j-q^{-l_j}t_{b_{i_1, i_2}}}{z_j-q^{l_j}t_{b_{i_1, i_2}}}
\right)
\prod_{i_2=s_{i_1}+1}^{\gamma_{i_1}}\frac{z_{k(i_1)}-q^{-l_{k(i_1)}}t_{b_{i_1, i_2}}}{z_{k(i_1)}-q^{l_{k(i_1)}}t_{b_{i_1, i_2}}}
\right\}
\\
\qquad{}\times
\left(\prod_{i_1<j_1}\frac{t_{b_{i_1, i_2}}-t_{b_{j_1, j_2}}}{t_{b_{i_1, i_2}}-q^{-2}t_{b_{j_1, j_2}}}\right),
\end{gather*}
where
\begin{gather*}
C_{(a) (\gamma)}
 = (-1)^{\sum\limits_{i=1}^r a_i}
\left(q^{\sum\limits_{s=1}^r\big(\sum\limits_{t=k(s)+1}^{n}l_t\big)(n_s-2 a_s)} \right)
\left(q^{(L-N)\big\{\sum\limits_{s=1}^{r}(n_s-2 a_s)\big\}}\right)
\nonumber
\\
\phantom{C_{(a) (\gamma)} =}{}
\times
\left(q^{-\sum\limits_{s=1}^{r}(n_s-2 a_s)\big(\sum\limits_{t=s+1}^r n_t\big)}\right)
\left(q^{2\sum\limits_{s=1}^{r}\sum\limits_{s<t}\gamma_t(a_s-\gamma_s)}\right)
\left(\prod_{1\le b \le N}q^{-1-\epsilon_b}\right)^{\sum\limits_{s=1}^{r}(a_s-\gamma_s)}.
\end{gather*}
Here we have used Lemma \ref{combi} $(i)$.

By $(ii)$ of  Lemma \ref{combi} we have
\begin{gather*}
\widehat{J}_{(\epsilon)(a)}^{(\nu)}
=
\sum_{j=1}^r\sum_{\substack{0\le \gamma_j \le a_j \\ 1 \le j \le r}} C_{(a) (\gamma)}
\sum_{\substack{1\le b_{i_1, i_2}\le N \\ 1\le i_1 \le r \\ 1 \le i_2 \le \gamma_{i_2}}}
\left(\prod_{i_1<j_1}\frac{t_{b_{i_1, i_2}}-t_{b_{j_1, j_2}}}{t_{b_{i_1, i_2}}-q^{-2}t_{b_{j_1, j_2}}}\right)
\\
\phantom{\widehat{J}_{(\epsilon)(a)}^{(\nu)}=}{} \times
\prod_{i_1=1}^{r}\left\{\sum_{s_{i_1}=0}^{\gamma_{i_1}}\sum_{\lambda_{i_1}=s_{i_1}}^{a_j-\gamma_{i_1}+s_{i_1}}
\sum_{m_{i_1}=0}^{n_{i_1}}(-1)^{m_{i_1}}
q^{-m_{i_1}(n_{i_1}-1)} q^{2l_{k(i_1)}(m_{i_1}-\lambda_{i_1})}\qbi{n_{i_1}}{m_{i_1}}
\right.
\\
\left.
\phantom{\widehat{J}_{(\epsilon)(a)}^{(\nu)}=}{}
\times
\left(\sum_{\substack{|A_{\mu, i_1}|=a_{i_1} \\ |A_{\mu, i_1}\cap M_{i_1}|=\lambda_{i_1}}}
\prod_{\substack{i_2<j_2 \\ i_1=j_1}}q^{-\mu_{i_1, i_2}}\right)
\left(
q^{\lambda_{i_1}\gamma_{i_1}+a_{i_1}\gamma_{i_1}-a_{i_1}s_{i_1}-\gamma_{i_1}^2}
\qbi{\lambda_{i_1}}{s_{i_1}}
\qbi{a_{i_1}-\lambda_{i_1}}{\gamma_{i_1}-s_{i_1}}
\right)
\right.
\\
\left.
\phantom{\widehat{J}_{(\epsilon)(a)}^{(\nu)}=}{}\times
\left\{\left(\prod_{i_2=1}^{\gamma_{i_1}}(1-q^{-1-\epsilon_{b_{i_1, i_2}}})
\prod_{b\ne b_{i_1, i_2}}\frac{t_{b_{i_1, i_2}}-q^{-1-\epsilon_b}t_b}{t_{b_{i_1, i_2}}-t_b}
\prod_{i_2<j_2}\frac{t_{b_{i_1, i_2}}-t_{b_{i_1, j_2}}}{t_{b_{i_1, i_2}}-q^{-2}t_{b_{i_1, j_2}}}
\right.
\right.
\right.
\\
\left.
\left.
\left.
\phantom{\widehat{J}_{(\epsilon)(a)}^{(\nu)}=}{}
\times
\prod_{j=1}^{k(i_1)-1}\frac{z_j-q^{-l_j}t_{b_{i_1, i_2}}}{z_j-q^{l_j}t_{b_{i_1, i_2}}}
\right)
\prod_{i_2=s_{i_1}+1}^{\gamma_{i_1}}\frac{z_{k(i_1)}-q^{-l_{k(i_1)}}t_{b_{i_1, i_2}}}{z_{k(i_1)}-q^{l_{k(i_1)}}t_{b_{i_1, i_2}}}
\right\}
\right\}
\\
\phantom{\widehat{J}_{(\epsilon)(a)}^{(\nu)}}{}=
\sum_{\substack{0\le \gamma_j \le a_j \\ 1 \le j \le r}} C_{(a) (\gamma)}
\sum_{\substack{1\le b_{i_1, i_2}\le N \\ 1\le i_1 \le r \\ 1 \le i_2 \le \gamma_{i_2}}}
\left(\prod_{i_1<j_1}\frac{t_{b_{i_1, i_2}}-t_{b_{j_1, j_2}}}{t_{b_{i_1, i_2}}-q^{-2}t_{b_{j_1, j_2}}}\right)
\\
\phantom{\widehat{J}_{(\epsilon)(a)}^{(\nu)}=}{}\times
\prod_{i_1=1}^{r}\left\{\sum_{s_{i_1}=0}^{\gamma_{i_1}}\sum_{\lambda_{i_1}=s_{i_1}}^{a_j-\gamma_{i_1}+s_{i_1}}
\sum_{m_{i_1}=0}^{n_{i_1}}(-1)^{m_{i_1}}
q^{-m_{i_1}(n_{i_1}-1)}
q^{2l_{k(i_1)}(m_{i_1}-\lambda_{i_1})}
\qbi{n_{i_1}}{m_{i_1}}
\right.
\\
\left.
\phantom{\widehat{J}_{(\epsilon)(a)}^{(\nu)}=}{}\times
\left(
q^{n_{i_1}\lambda_{i_1}+a_{i_1}n_{i_1}-a_{i_1}m_{i_1}-a_{i_1}-\frac{n_{i_1}(n_{i_1}-1)}{2}}
\qbi{m_{i_1}}{\lambda_{i_1}}
\qbi{n_{i_1}-m_{i_1}}{a_{i_1}-\lambda_{i_1}}
\right)
\right.
\\
\left.
\phantom{\widehat{J}_{(\epsilon)(a)}^{(\nu)}=}{}\times
\left(
q^{\lambda_{i_1}\gamma_{i_1}+a_{i_1}\gamma_{i_1}-a_{i_1}s_{i_1}-\gamma_{i_1}^2}
\qbi{\lambda_{i_1}}{s_{i_1}}
\qbi{a_{i_1}-\lambda_{i_1}}{\gamma_{i_1}-s_{i_1}}
\right)
\right.
\\
\left.
\phantom{\widehat{J}_{(\epsilon)(a)}^{(\nu)}=}{}\times
\left\{\prod_{i_2=1}^{\gamma_{i_1}}\left((1-q^{-1-\epsilon_{b_{i_1, i_2}}})
\prod_{b\ne b_{i_1, i_2}}\frac{t_{b_{i_1, i_2}}-q^{-1-\epsilon_b}t_b}{t_{b_{i_1, i_2}}-t_b}
\prod_{i_2<j_2}\frac{t_{b_{i_1, i_2}}-t_{b_{i_1, j_2}}}{t_{b_{i_1, i_2}}-q^{-2}t_{b_{i_1, j_2}}}
\right.
\right.
\right.
\\
\left.
\left.
\left.
\phantom{\widehat{J}_{(\epsilon)(a)}^{(\nu)}=}{}\times
\prod_{j=1}^{k(i_1)-1}\frac{z_j-q^{-l_j}t_{b_{i_1, i_2}}}{z_j-q^{l_j}t_{b_{i_1, i_2}}}
\right)
\prod_{i_2=s_{i_1}+1}^{\gamma_{i_1}}\frac{z_{k(i_1)}-q^{-l_{k(i_1)}}t_{b_{i_1, i_2}}}{z_{k(i_1)}-q^{l_{k(i_1)}}t_{b_{i_1, i_2}}}
\right\}
\right\}.
\end{gather*}
It is easy to show
\begin{gather*}
 \sum_{\lambda=s}^{a-\gamma+s}
\sum_{m=0}^{n}(-1)^{m}
q^{-m(n-1)}
\qbi{n}{m}
q^{2l(m-\lambda)}
\left(q^{n \lambda +a n-a m-a-\frac{n(n-1)}{2}}
\qbi{m}{\lambda}
\qbi{n-m}{a-\lambda}
\right)
\\
\qquad{}\times
\left(
q^{\lambda \gamma+a \gamma-a s-\gamma^2}
\qbi{\lambda}{s}
\qbi{a-\lambda}{\gamma-s}
\right)
\\
\qquad{}=
(-1)^s
q^{a(n-s-1)+s}
q^{-\frac{n(n-1)}{2}}
\frac{[n]!}{[s]![a-s]!}
\sum_{i=0}^{n-a}
(-1)^{i}q^{i(2l-n-a+1)}
\frac{1}{[i]![n-a-i]!}\delta_{a, \gamma},
\end{gather*}
for $0\le s \le \gamma \le a\le n$.

Hence
\begin{gather*}
\widehat{J}_{(\epsilon)(a)}^{(\nu)}
=(-1)^{\sum\limits_{i=1}^r a_i}
\left(q^{\sum\limits_{s=1}^r\big(\sum\limits_{t=k(s)+1}^{n}l_t\big)(n_s-2 a_s)} \right)
\left(q^{(L-N)\big\{\sum\limits_{s=1}^{r}(n_s-2 a_s)\big\}}\right)
\\
\phantom{\widehat{J}_{(\epsilon)(a)}^{(\nu)}=}{} \times
\left(q^{-\sum\limits_{s=1}^{r}(n_s-2 a_s)\big(\sum\limits_{t=s+1}^r n_t\big)}\right)
\sum_{\substack{1\le b_{i_1, i_2}\le N \\ 1\le i_1 \le r \\ 1\le i_2 \le a_{i_1}}}
\left(\prod_{i_1<j_1}\frac{t_{b_{i_1, i_2}}-t_{b_{j_1, j_2}}}{t_{b_{i_1, i_2}}-q^{-2}t_{b_{j_1, j_2}}}\right)
\\
\phantom{\widehat{J}_{(\epsilon)(a)}^{(\nu)}=}{}\times
\prod_{i_1=1}^{r}\left\{\sum_{s_{i_1}=0}^{a_{i_1}}
(-1)^{s_{i_1}}
q^{a_{i_1}(n_{i_1}-s_{i_1}-1)+s_{i_1}}
q^{-\frac{n_{i_1}(n_{i_1}-1)}{2}}
\frac{[n_{i_1}]!}{[s_{i_1}]![a_{i_1}-s_{i_1}]!}
\right.
\\
\left.\phantom{\widehat{J}_{(\epsilon)(a)}^{(\nu)}=}{}\times
\sum_{i=0}^{n_{i_1}-a_{i_1}}
(-1)^{i}q^{i(2l_{k(i_1)}-n_{i_1}-a_{i_1}+1)}
\frac{1}{[i]![n_{i_1}-a_{i_1}-i]!}
\right.
\\
\left.
\phantom{\widehat{J}_{(\epsilon)(a)}^{(\nu)}=}{}\times
\left\{\prod_{i_2=1}^{a_{i_1}}\left((1-q^{-1-\epsilon_{b_{i_1, i_2}}})
\prod_{b\ne b_{i_1, i_2}}\frac{t_{b_{i_1, i_2}}-q^{-1-\epsilon_b}t_b}{t_{b_{i_1, i_2}}-t_b}
\prod_{i_2<j_2}\frac{t_{b_{i_1, i_2}}-t_{b_{i_1, j_2}}}{t_{b_{i_1, i_2}}-q^{-2}t_{b_{i_1, j_2}}}
\right.
\right.
\right.
\\
\left.
\left.
\left.
\phantom{\widehat{J}_{(\epsilon)(a)}^{(\nu)}=}{}\times
\prod_{j=1}^{k(i_1)-1}\frac{z_j-q^{-l_j}t_{b_{i_1, i_2}}}{z_j-q^{l_j}t_{b_{i_1, i_2}}}
\right)
\prod_{i_2=s_{i_1}+1}^{a_{i_1}}\frac{z_{k(i_1)}-q^{-l_{k(i_1)}}t_{b_{i_1, i_2}}}{z_{k(i_1)}-q^{l_{k(i_1)}}t_{b_{i_1, i_2}}}
\right\}
\right\}.\tag*{\qed}
\end{gather*}\renewcommand{\qed}{}
\end{proof}

\begin{lemma}
If $a_i \ne n_i $ for some $i$,
\begin{gather*}
\sum_{\substack{\epsilon_i=\pm \\ 1\le i \le N}}\left(\prod_{j=1}^N \epsilon_j \right)
\left(\prod_{a<b}\frac{q^{\epsilon_b}t_b-q^{\epsilon_a}t_a}{t_b-q^{-2}t_a}\right)\widehat{J}_{(\epsilon)(a)}^{(\nu)}
=0.
\end{gather*}
\end{lemma}

\begin{proof} It is enough to show the following equation.
For $1\le b_{i_1, i_2}\le N$ $(1\le i_1 \le r$, $1\le i_2 \le a_{i_1})$,
$b_{i_1, i_2}\ne b_{j_1, j_2} $ $((i_1, i_2)\ne(j_1, j_2))$,
\begin{gather}
\sum_{\substack{\epsilon_i=\pm \\ 1\le i \le N}}\left(\prod_{i=1}^N \epsilon_{i}\right)
\prod_{a<b}(q^{\epsilon_b}t_b-q^{\epsilon_a}t_a)
\prod_{\substack{1\le i_1 \le r \\ 1\le i_2 \le a_{i_1}}}
\left\{(1-q^{-1-\epsilon_{b_{i_1, i_2}}})
\prod_{b\ne b_{i_1, i_2}}(t_{b_{i_1, i_2}}-q^{-1-\epsilon_b}t_b)
\right\}
\nonumber
\\
\qquad =(1-q^{-2})^{N} q^{\frac{N(N-1)}{2}}\left(\prod_{s=1}^{r}\delta_{a_s, n_s}\right)
\left\{\prod_{a<b} (t_b-t_a)\right\}
\left\{\prod_{b\ne b_{i_1, i_2}}(t_{b_{i_1, i_2}}-q^{-2}t_b)\right\}.
\label{mod}
\end{gather}
For a set $\{b_{1, 1}, \dots, b_{r, a_r}\}=\{b_1, \dots, b_{\alpha}\}$,
let$\{c_1, \dots ,c_{N-\alpha}\}$ be def\/ined by
\begin{gather*}
\{b_1, \dots, b_{\alpha}\}\sqcup \{c_1, \dots ,c_{N-\alpha}\}=\{1, \dots, N \} ,
\end{gather*}
where $\alpha=\sum\limits_{i=1}^r a_i$.

Then the left hand side of (\ref{mod}) is equal to
\begin{gather*}
(1-q^{-2})^{\alpha}\left(\prod_{1\le i \le \alpha}\delta_{\epsilon_{b_i}, +} \right)
\left\{\prod_{i<j}q (t_{b_j}- t_{b_i})\right\}
\left\{\prod_{\substack{1\le i, j \le \alpha \\ i\ne j}}
(t_{b_{i}}-q^{-2}t_{b_j})\right\}
\\
\qquad{}
\times
\left\{\prod_{b_i<c_j}(-q)\right\}
\left\{\prod_{c_i<b_j}q \right\}
\sum_{\substack{\epsilon_{c_i}=\pm \\ 1\le i \le N-\alpha}}\left(\prod_{i=1}^{N-\alpha} \epsilon_{c_i}\right)
\left\{\prod_{i<j}(q^{\epsilon_{c_j}}t_{c_j}-q^{\epsilon_{c_i}}t_{c_i})\right\}
\\
\qquad{}
\times
\left\{\prod_{1\le i \le \alpha}\prod_{1\le j \le N-\alpha}(t_{b_i}-q^{\epsilon_{c_j}-1}t_{c_j})\right\}
\left\{\prod_{1\le i \le \alpha}\prod_{1\le j \le N-\alpha}(t_{b_{i}}-q^{-1-\epsilon_{c_j}}t_{c_j})\right\}
.
\end{gather*}
Using
\begin{gather*}
(t_{b_j}-q^{\epsilon_{c_i}-1}t_{c_i})\big(t_{b_{i}}-q^{-1-\epsilon_{c_j}}t_{c_j}\big)
=(t_{b_j}-t_{c_i})\big(t_{b_{i}}-q^{-2}t_{c_j}\big),
\end{gather*}
we have
\begin{gather*}
\sum_{\epsilon}\left(\prod_{i=1}^N \epsilon_{i}\right)\!
\left\{\prod_{a<b}\big(q^{\epsilon_b}t_b-q^{\epsilon_a}t_a\big)\right\}\!
\left\{\prod_{\substack{1\le i_1 \le r \\ 1\le i_2 \le a_{i_1}}}
\big(1-q^{-1-\epsilon_{b_{i_1, i_2}}}\big)\!\right\}\!
\left\{\prod_{b\ne b_{i_1, i_2}}\big(t_{b_{i_1, i_2}}-q^{-1-\epsilon_b}t_b\big)\!\right\}\!
\\
\qquad {} =(1-q^{-2})^{\alpha}\left(\prod_{1\le i \le \alpha}\delta_{\epsilon_{b_i}, +} \right)
\!\left\{\prod_{i<j}q (t_{b_j}- t_{b_i})\!\right\}\!
\left\{\prod_{\substack{1\le i, j \le \alpha \\ i\ne j}}
(t_{b_{i}}-q^{-2}t_{b_j})\!\right\}\!
\\
\qquad{}
\times
\left\{\prod_{1\le i \le \alpha}\prod_{1\le j \le N-\alpha}(t_{b_i}-t_{c_j})
\big(t_{b_{i}}-q^{-2}t_{c_j}\big)\right\}
\\
\qquad{}
\times
\left\{\prod_{b_i<c_j}(-q)\right\}
\left\{\prod_{c_i<b_j}q \right\}
\sum_{\genfrac{}{}{0pt}{}{\epsilon_{c_i}=\pm}{1\le i \le N-\alpha}}\left(\prod_{i=1}^{N-\alpha} \epsilon_{c_i}\right)
\left\{\prod_{i<j}\big(q^{\epsilon_{c_j}}t_{c_j}-q^{\epsilon_{c_i}}t_{c_i}\big)\right\}
.
\end{gather*}
Let $\alpha\ne N$ and
${\boldsymbol{a}}_i(\epsilon)=^{t}(1, q^{\epsilon}t_i, (q^{\epsilon}t_i)^2, \dots, (q^{\epsilon}t_i)^{N-\alpha-1})$.
Then
\begin{gather}
\sum_{\substack{\epsilon_{i}=\pm \\ 1\le i \le N-\alpha}}\left(\prod_{i=1}^{N-\alpha} \epsilon_{i}\right)
\prod_{i<j}(q^{\epsilon_{j}}t_{j}-q^{\epsilon_{i}}t_{i})\nonumber\\
\qquad{}
=\sum_{\substack{\epsilon_{i}=\pm\\ 1\le i \le N-\alpha}}\left(\prod_{i=1}^{N-\alpha} \epsilon_{i}\right)
\det ({\boldsymbol{a}}_1(\epsilon_1), {\boldsymbol{a}}_2(\epsilon_2), \dots, {\boldsymbol{a}}_{N-\alpha}(\epsilon_{N-\alpha})).
\label{0}
\end{gather}
Since
\begin{gather*}
\sum_{\epsilon_i=\pm}\epsilon_i{\boldsymbol{a}}_i(\epsilon)=
^{t}(0, (q-q^{-1})t_i, \dots, (q^{N-\alpha-1}-q^{-(N-\alpha-1)})t_i^{N-\alpha-1}),
\end{gather*}
the right hand side of (\ref{0}) is equal to $0$.
\end{proof}

If $a_i=n_i$ for all $i$, then
\begin{gather}
\sum_{\substack{\epsilon_i=\pm \\ 1\le i \le N}}\left(\prod_{i=1}^N \epsilon_i\right)
\left(\prod_{1\le a< b\le N}\frac{q^{\epsilon_b} t_b-q^{\epsilon_a} t_a}{t_b-q^{-2} t_a}\right)\hat{J}_{(\epsilon)(a)}^{(\nu)}\nonumber\\
\qquad{}=
C_1\left(\prod_{a<b} \frac{t_b-t_a}{t_b-q^{-2}t_a}\right)
\sum_{\substack{\Gamma_1\sqcup \dots \sqcup \Gamma_r=\{1, \dots, N\}\\ |\Gamma_s|=n_s\,\,(s=1, \dots, r)}}
\sum_{\substack{b_{i_1, i_2}\in \Gamma_{i_1} \\ 1\le i_1 \le r \\ 1\le i_2 \le n_{i_1}}}
\left(\prod_{i_1> j_1}\frac{t_{b_{i_1, i_2}}-q^{-2}t_{b_{j_1, j_2}}}{t_{b_{i_1, i_2}}-t_{b_{j_1, j_2}}}\right)
\nonumber
\\
\qquad{}
\times
\prod_{i_1=1}^{r}\left\{
\sum_{s_{i_1}=0}^{n_{i_1}}
(-1)^{s_{i_1}}
q^{-(n_{i_1}-1)s_{i_1}}
\qbi{n_{i_1}}{s_{i_1}}
\prod_{i_2=1}^{s_{i_1}}(z_{k(i_1)}-q^{l_{k(i_1)}}t_{b_{i_1, i_2}})
\right.
\nonumber
\\
\left.\qquad{}
\times
\prod_{i_2=s_{i_1}+1}^{n_{i_1}}(z_{k(i_1)}-q^{-l_{k(i_1)}}t_{b_{i_1, i_2}})
\prod_{i_2>j_2}\frac{t_{b_{i_1, i_2}}-q^{-2}t_{b_{i_1, j_2}}}{t_{b_{i_1, i_2}}-t_{b_{i_1, j_2}}}
\right.
\nonumber
\\
\left.\qquad{}
\times
\prod_{i_2=1}^{n_{i_1}}\left(
\frac{1}{z_{k(i_1)}-q^{l_{k(i_1)}}t_{b_{i_1, i_2}}}
\prod_{j=1}^{k(i_1)-1}\frac{z_j-q^{-l_j}t_{b_{i_1, i_2}}}{z_j-q^{l_j}t_{b_{i_1, i_2}}}
\right)
\right\},
\label{a=n}
\end{gather}
where
\begin{gather*}
C_1=(-1)^N\big(1-q^{-2}\big)^{N} q^{N^2-LN} q^{\frac{N(N-1)}{2}} q^{\sum\limits_{i=1}^r \frac{n_i(n_i-1)}{2}}\\
\phantom{C_1=}{}\times
\left\{q^{-\sum\limits_{s=1}^r\big(\sum\limits_{t=k(s)+1}^{n}l_t\big) n_s} \right\}
\left\{q^{\sum\limits_{s=1}^{r}\big(\sum\limits_{t=s+1}^r n_t\big)n_s}\right\}.
\end{gather*}
By Lemma {\ref{z}} the right hand side of (\ref{a=n}) becomes
\begin{gather*}
 C_1\prod_{s=1}^r\left\{(-1)^{n_s}[n_s]! q^{-l_{k(s)}n_s} q^{-\frac{n_s(n_s-1)}{2}}
\left\{\prod_{i=0}^{n_s-1}\big(1-q^{2(l_{k(s)}-i)}\big)\right\}\right\}
\\
\qquad{}
\times
\left(\prod_{a<b} \frac{t_b-t_a}{t_b-q^{-2}t_a}\right)
\sum_{\genfrac{}{}{0pt}{}{\Gamma_1\sqcup \dots \sqcup \Gamma_r=\{1, \dots, N\}}{|\Gamma_s|=n_s\,\,(s=1, \dots, r)}}
\left(\prod_{\genfrac{}{}{0pt}{}{1\le i<j \le r}{a\in \Gamma_i,  b\in \Gamma_j }}\frac{t_{b}-q^{-2}t_{a}}{t_{b}-t_{a}}\right)
\\
\qquad{}
\times
\prod_{s=1}^{r}
\prod_{b\in \Gamma_s}
\left(
\frac{t_b}{z_{k(s)}-q^{l_{k(s)}}t_{b}}
\prod_{i=1}^{k(s)-1}\frac{z_i-q^{-l_i}t_{b}}{z_i-q^{l_i}t_{b}}
\right)
.
\end{gather*}
This completes the proof of Proposition $\ref{mainprop}$.\hfill \qed

\appendix
\section[The representation $V^{(l)}_z$]{The representation $\boldsymbol{V^{(l)}_z}$}\label{appendixA}

 Let $q^{h_i}$, $e_i$, $f_i$ $(i=0, 1)$ and $q^d$ be the generators of~$U_q(\widehat{sl_2})$. (See~\cite{idzumi} for more details.)
The actions of the generators of $U_q(\widehat{sl_2})$ on $V^{(l)}_z$ are given as follows.

For $0\le i \le l$ and $n \in {\mathbb{Z}}$,
\begin{alignat*}{3}
& e_0 v_j^{(l)}\otimes z^n=[l-i]v_{i+1}^{(l)}\otimes z^{n+1},\qquad
&&
e_1 v_j^{(l)}\otimes z^n=[i]v_{i-1}^{(l)}\otimes z^{n},&
\\
& f_0 v_j^{(l)}\otimes z^n=[i]v_{i-1}^{(l)}\otimes z^{n-1},\qquad &&
f_1 v_j^{(l)}\otimes z^n=[l-i]v_{i+1}^{(l)}\otimes z^n,&
\\
& q^{h_0} v_i^{(l)}\otimes z^n=q^{-(l-2i)}v_i^{(l)}\otimes z^n,\qquad &&
q^{h_1} v_i^{(l)}\otimes z^n=q^{l-2i}v_i^{(l)}\otimes z^n,&
\\
& q^d v_j^{(l)}\otimes z^n=q^n v_j^{(l)}\otimes z^n. &&&
\end{alignat*}

\section[$R$-matrix]{$\boldsymbol{R}$-matrix}\label{appendixB}

We give examples of explicit forms of $R$-matrix in the case of $l_1=1$ or $l_2=1$. They are taken from~\cite{idzumi}.
If we write
\begin{gather*}
R_{1, l_2}(z) \big(v_{\epsilon}^{(1)}\otimes v_j^{(l_2)}\big)=\sum_{\epsilon'=0, 1}
 v_{\epsilon'}^{(1)}\otimes r_{\epsilon' \epsilon}^{1 l_2}(z)v_j^{(l_2)},\\
R_{l_1, 1}(z) \big(v_j^{(l_1)}\otimes v_{\epsilon}^{(1)}\big)=\sum_{\epsilon'=0, 1}
r_{\epsilon' \epsilon}^{l_1 1}(z)v_j^{(l_1)}\otimes v_{\epsilon'}^{(1)} ,
\end{gather*}
then we have
\begin{gather*}
\left(
\begin{array}{cc}
r_{0 0}^{1 l_2}(z)& r_{0 1}^{1 l_2}(z)\\
r_{1 0}^{1 l_2}(z)& r_{1 1}^{1 l_2}(z)
\end{array}
\right)
 = \frac{1}{q^{1+l_2/2}-z^{-1}q^{-l_2/2}}\left(
\begin{array}{cc}
q^{1+h/2}-z^{-1}q^{-h/2} & (q-q^{-1})z^{-1}fq^{h/2}\\
(q-q^{-1})eq^{-h/2}& q^{1-h/2}-z^{-1}q^{h/2}
\end{array}
\right)
,
\\
\left(
\begin{array}{cc}
r_{0 0}^{l_1 1}(z)& r_{0 1}^{l_1 1}(z)\\
r_{1 0}^{l_1 1}(z)& r_{1 1}^{l_1 1}(z)
\end{array}
\right)
=\frac{1}{zq^{l_1/2}-q^{-1-l_1/2}}\left(
\begin{array}{cc}
z q^{h/2}-q^{-1-h/2} & (q-q^{-1}) z q^{h/2}f\\
(q-q^{-1}) q^{-h/2} e& z q^{-h/2}-q^{-1+h/2}
\end{array}
\right)
,
\end{gather*}
$h=h_1$, $e=e_1$ and $f=f_1$.

\section[Free field representations]{Free f\/ield representations}\label{appendixC}

The following formulae are given in \cite{kato}. For $x=a, b, c$ let
\begin{gather*}
x(L; M, N|z: \alpha) = -\sum_{n\ne 0}\frac{[Ln]x_n}{[Mn][Nn]}z^{-n}q^{|n|\alpha}+\frac{L\tilde{x}_0}{MN}\log z
+\frac{L}{MN}Q_x,
\\
x(N|z: \alpha)=x(L; L, N|z: \alpha)
=-\sum_{n\ne 0}\frac{x_n}{[Nn]}z^{-n}q^{|n|\alpha}+\frac{\tilde{x}_0}{N}\log z
+\frac{1}{N}Q_x.
\end{gather*}
The normal ordering is def\/ined by specifying $N_+$, $\tilde{a}_0$, $\tilde{b}_0$, $\tilde{c}_0$ as annihilation operators,
$N_-$, $Q_a$, $Q_b$, $Q_c$ as creation operators.

Def\/ine operators
\begin{gather*}
 J^{-}(z): \ F_{r, s}\to F_{r, s+1},
\qquad
 S(z):\ F_{r, s}\to F_{r-2, s-1},
\qquad
 \phi_m^{(l)}(z): \ F_{r, s}\to F_{r+l, s+l-m},
\end{gather*}
by
\begin{gather*}
J^{-}(z) = \frac{1}{(q-q^{-1})z}(J^{-}_+(z)-J^{-}_-(z)),
\\
J^{-}_{\mu}(z)=:\exp\left(a^{(\mu)}\left(q^{-2}z;-\frac{k+2}{2}\right)+b\big(2|q^{(\mu-1)(k+2)}z; -1\big)
+c\big(2|q^{(\mu-1)(k+1)-1}z;0\big)\right):,
\\
a^{(\mu)}\left(q^{-2}z; -\frac{k+2}{2}\right)
 = \mu\left\{\big(q-q^{-1}\big)\sum_{n=1}^{\infty}a_{\mu n}z^{-\mu n}q^{(2\mu-\frac{k+2}{2})n}+\tilde{a}_0 \log q
\right\},
\\
S(z)=\frac{-1}{(q-q^{-1})z}(S_+(z)-S_-(z)),
\nonumber
\\
S_{\epsilon}(z)=:\exp\left(-a\left(k+2|q^{-2}z; -\frac{k+2}{2}\right)-b\big(2|q^{-k-2}z;-1\big)-c\big(2|q^{-k-2+\epsilon}z; 0\big)\right):,
\\
\phi^{(l)}_l(z) = :\exp\left(a\left(l; 2, k+2|q^k z; \frac{k+2}{2}\right)\right):,
\\
\phi^{(l)}_{l-r}(z) = \frac{1}{[r]!}\oint\left(\prod_{j=1}^r\frac{d u_j}{2\pi i}\right) \left[\dots\left[\left[\phi^{(l)}_l(z), J^-(u_1)\right]_{q^l}, J^-(u_2)\right]_{q^{l-2}},\dots ,
 J^{-}(u_r)\right]_{q^{l-2r+2}},
\end{gather*}
where
\begin{gather*}
[r]!=\prod_{i=1}^{r}[i],\qquad [X, Y]_q=XY-qYX,
\end{gather*}
and the integral in $\phi^{(l)}_{l-r}(z)$ signif\/ies to take the coef\/f\/icient of $(u_1\cdots u_r)^{-1}$.

\section{List of OPE's}\label{appendixD}

The following formulae are given in \cite{kato}
\begin{gather*}
 \phi_{l_1}(z_1) \phi_{l_2}(z_2)=(q^k z_1)^{\frac{l_1 l_2}{2(k+2)}}
\frac{\left(q^{l_1+l_2+2k+6}\frac{ z_2}{z_1} ; q^4 , q^{2(k+2)}\right)_{\infty}
\left(q^{-l_1-l_2+2k+6}\frac{ z_2}{z_1} ; q^4 , q^{2(k+2)}\right)_{\infty}}
{\left(q^{l_1-l_2+2k+6} \frac{z_2}{z_1} ; q^4 , q^{2(k+2)}\right)_{\infty}
\left(q^{-l_1+l_2+2k+6} \frac{z_2}{z_1} ; q^4 , q^{2(k+2)}\right)_{\infty}}
\\
\phantom{\phi_{l_1}(z_1) \phi_{l_2}(z_2)=}{} \times
:\phi_{l_1}(z_1) \phi_{l_2}(z_2):,
\qquad |q^{-l_1-l_2+2k+6}z_2|<|z_1|,
\\
\phi_l(z) J^{-}_{\mu}(u)=\frac{z-q^{\mu l-k-2}u}{z-q^{l-k-2}u}:\phi_l(z) J^{-}_{\mu}(u):,\qquad
|q^{-l-k-2}u|<|z|,
\\
J^{-}_{\mu}(u) \phi_l(z)=q^{\mu l}\frac{u-q^{-\mu l+k+2}z}{u-q^{l+k+2}z}:\phi_l(z) J^{-}_{\mu}(u):,\qquad
|q^{-l+k+2}u|<|z|,
\\
 \phi_l(z) S_{\epsilon}(t)=\frac{\left(q^l\frac{t}{z}; p\right)_{\infty}}{\left(q^{-l}\frac{t}{z}; p\right)_{\infty}}
(q^k z)^{-\frac{l}{k+2}} :\phi_l(z) S_{\epsilon}(t):,\qquad |z|>|q^{-l} t|,
\\
J^{-}_{\mu}(u) S_{\epsilon}(t)=q^{-\mu}\frac{u-q^{-\mu(k+1)-\epsilon} t}{u-q^{-\mu(k+2)} t}
:J^{-}_{\mu}(u) S_{\epsilon}(t):,\qquad |u|>|q^{-k-2} t|,
\\
J^{-}_{\mu_1}(u_1) J^{-}_{\mu_2}(u_2)=\frac{q^{-\mu_1} u_1-q^{-\mu_2}u_2}{u_1-q^{-2}u_2} :J^{-}_{\mu_1}(u_1) J^{-}_{\mu_2}(u_2):,\qquad |u_1|>|q^{-2} u_2|,
\\
S_{\epsilon_1}(t_1)S_{\epsilon_2}(t_2)=(q^{-2} t_1)^{\frac{2}{k+2}} \frac{q^{\epsilon_1} t_1-q^{\epsilon_2} t_2}{t_1-q^{-2} t_2}
\frac{\left(q^{-2}\frac{t_2}{t_1}; p\right)_{\infty}}{\left(q^{2}\frac{t_2}{t_1}; p\right)_{\infty}}
:S_{\epsilon_1}(t_1)S_{\epsilon_2}(t_2):,\qquad |t_1|>|q^{-2} t_2|.
\end{gather*}

\subsection*{Acknowledgements}

After completing the paper we know that similar results to the present
paper are obtained by H.~Awata, S.~Odake and J.~Shiraishi \cite{shiraishi2}.
I would like to thank all of them for kind correspondences and permitting me
to write the results in a single authored paper.
I am also grateful to Professor H.~Konno for useful comments.
I would like to thank Professor A.~Nakayashiki for constant encouragements.

\pdfbookmark[1]{References}{ref}

\LastPageEnding

\end{document}